\documentclass{article}
\usepackage{curves}
\usepackage{amsmath}
\usepackage{amsfonts}
\usepackage{amsthm}
\usepackage{graphicx}
\usepackage{float}
\usepackage{url}
\usepackage{framed}
\usepackage[top= 1.25in, bottom=1.25in, left=1.15in, right=1.15in]{geometry}
\usepackage{tikz}
\usepackage{sans}
\usetikzlibrary{arrows,snakes,backgrounds,decorations.pathreplacing,decorations.markings}

\title{{Classifying (almost)-Belyi maps with Five Exceptional Points}}
\newtheorem{theorem}{Theorem}[section]
\newtheorem{lemma}[theorem]{\textbf{Lemma}}
\newtheorem{proposition}[theorem]{\textbf{Proposition}}

\newtheorem{remark}[theorem]{\textbf{Remark}}
\newtheorem{defn}[theorem]{\textbf{Definition}}
\newtheorem{exmp}[theorem]{\textbf{Example}}

\newtheorem{Conjecture}{Conjecture}
\newcommand{\eindebewijs}{\hfill$\Box$\par\medskip}

\begin{document}
\author{Mark van Hoeij\thanks{Supported by NSF grants 1319547 and 1618657.} \\
{Department of Mathematics} \\
{Florida State University} \\
{Tallahassee, FL 32306, USA}
\and
Vijay Jung Kunwar$^*$ \\
{Department of Mathematics} \\
{Albany State University} \\
{Albany, GA 31705, USA}
}

\maketitle
\begin{abstract}
\label{abstract}
We classify all rational functions $f: \mathbb{P}^1\rightarrow \mathbb{P}^1$ whose branching pattern above $0,1,\infty$ satisfy
a certain regularity condition with precisely $d=5$ exceptions.
This work is motivated by solving second order linear differential equations, with $d=5$ true singularities, in terms of hypergeometric functions.
A similar problem was solved for $d=4$ in \cite{vidunas bound}. 
\end{abstract}

\section{Introduction}
\label{newIntro}
Our main goal in this paper is to tabulate all rational functions $f \in \mathbb{C}(x)$ whose
branching patterns satisfy a regularity condition defined in Section \ref{introduction}, and to prove completeness of the table.
This condition comes from solving differential equations with at most $d=5$ true singularities in terms of hypergeometric functions. The cases $d=3$ resp. $d=4$
were previously studied in \cite{vidunas-d3} resp. \cite{hoeijheun, vidunas bound}.
The functions $f$ in our tables are either Belyi maps or almost-Belyi maps:

\begin{defn}
A holomorphic map $f$ from a compact Riemann surface $C$ to the Riemann sphere \linebreak $\mathbb{P}^1 = \mathbb{C} \bigcup \{\infty\}$
is a {\em Belyi map} if its branched set is in $\{0, 1, \infty\}$, i.e. $f$ is unramified outside of $\{0, 1, \infty\}$.
\end{defn}
The pre-image $\subset C$ of the closed interval $[0,1] \subset \mathbb{P}^1$
under a Belyi map $f$ gives a bi-colored oriented graph, called dessin d'enfant.
Up to equivalence, there is a 1-1 correspondence between dessins d'enfants, 3-constellations, and Belyi maps, see Section \ref{section3} for details.
In our application $C = \mathbb{P}^1$ since our $f$'s are rational functions.

 \emph{Almost-Belyi maps} \cite{almost-belyi} are rational maps with only 1 or 2 simple branch points outside of $\{0, 1, \infty\}$. We denote these as $\text{Belyi}^{(1)}$ resp. $\text{Belyi}^{(2)}$ maps if they have 1 resp. 2 simple branch points outside $\{0, 1, \infty\}$.
Those branch points in $\mathbb{P}^1 \setminus \{0, 1, \infty\}$
are free to move, so while Belyi-maps are classified by discrete objects (dessins d'enfants, or 3-constellations),
$\text{Belyi}^{(1)}$ resp. $\text{Belyi}^{(2)}$ maps naturally occur in 1 resp. 2 dimensional families.

We expect our tables to be helpful in other contexts as well; 
Belyi maps have wide range of application in the fields of algebra, geometry, and combinatorics.
They are used to prove Davenport-Stothers-Zannier bound \cite{Belyi maps Applications}.
Shabat polynomials are the special case of Belyi maps with only one pole at infinity.
A \emph{dynamical} Belyi map $f: \mathbb{P}^1 \rightarrow \mathbb{P}^1$, which sends $\{0, 1, \infty \}$ to
$\{0, 1, \infty \}$, is used to construct the situation of a Julia set with ``complete chaos" \cite{constellation}.
Another application is the classification of elliptic fibrations of K3-surfaces \cite{beukers}. Some $\text{Belyi}^{(1)}$ maps correspond to solutions of isomonodromic systems of Fuchsian equations with 4 (+1 apparent) singularities, and hence give algebraic solutions of the Painlev\'{e} VI equation by Jimbo-Miwa correspondence \cite{Jimbo-Miwa}. Such $\text{Belyi}^{(1)}$ maps are studied as deformations of dessins d'infants in \cite{Kit1}.
The motivation for this paper is solving differential equations in terms of hypergeometric functions \cite{hoeijheun,hypergeomdeg3,vidunas bound,VJKThesis, RamifiedCoverings}
but we expect that our tables will be useful for other applications as well. Indeed, some entries of the tables have already appeared in prior applications, see Section~\ref{introduction}
for more.

\subsection{Motivation, solving differential equations}
\label{motivationIntro}

\begin{Conjecture} \label{Conj1}
Let $L$ be a second order linear differential equation $a_2 y'' + a_1 y' + a_0 y = 0$ with coefficients $a_i \in \mathbb{C}[x]$.
If $L$ is {\em regular singular}\footnote{We focus on the regular singular case
because for irregular singular equations of order 2, a complete algorithm to find
all $\{$Airy, Bessel, Kummer, Whittaker$\}$-type solutions was given
in~\cite{QThesis,debeerst}. The regular-singular assumption can be replaced by the assumption that $y$
in Conjecture~\ref{Conj1} has a non-zero radius of convergence.} (i.e. Fuchsian) and
has, among its solutions, a non-zero power series solutions with integer coefficients, ($y \in \mathbb{Z}[[x]] - \{0\}$), then one of these cases holds:
\begin{itemize}
\item $y$ is an {\em algebraic} function, or 
\item $y$ can be written as 
$y = 
r_0S(f)+r_1 (S(f))'$,  
where $S(f)=\mbox{}_2F_1(a, b; c \,|\,f)$.
Here $f,r_0,r_1$ are algebraic functions,
$a,b \in \mathbb{Q}$ and $c$ is a positive integer\footnote{This condition implies at least one logarithmic singularity. Because of the conjecture we focus on differential equations with
at least one logarithmic singularity, however, the same table can also be used for more general ``parametric'' cases \cite{vidunas bound}.}.
\end{itemize}
\end{Conjecture}

If $y$ is algebraic, it can be found with Kovacic' algorithm \cite{kovacic}, so we
are mainly interested in $\mbox{}_2F_1$-type solutions.
We only treat rational $f$'s in this paper, and plan to treat algebraic $f$'s later by exploiting their relation to the modular curve $X_0(N)$.
For $d=3$ singularities $f$ is a M\"{o}bius transformation, and for $d=4$ singularities (Heun's case) $f$'s are classified in \cite{vidunas bound,hoeijheun}. So we treat $d=5$ in this paper.

Recent algorithms for finding closed form solutions (solutions expressible in terms of well studied special functions)
are given in \cite{hoeijheun,QThesis,hypergeomdeg3,VJKThesis,debeerst,tfang}. Although random differential equations are unlikely to have closed form solutions, the conjecture says that the second order differential equations that are of most interest to combinatorics\footnote{Equations with a convergent integer power series solution ({\em globally nilpotent} differential equations).}
should have closed form solutions. 
We tested this on numerous differential equations obtained from the oeis.org (the Online Encyclopedia of Integer Sequences).
All turned out to be $\mbox{}_2F_1$-solvable with parameters that can be related 
to a
triplet $(k,\ell,m)$ (see the notation in Sections \ref{introduction} and \ref{motivation}) in Diagram~(1) in Takeuchi's classification \cite[Section 4]{takeuchi} of arithmetic triangle groups:

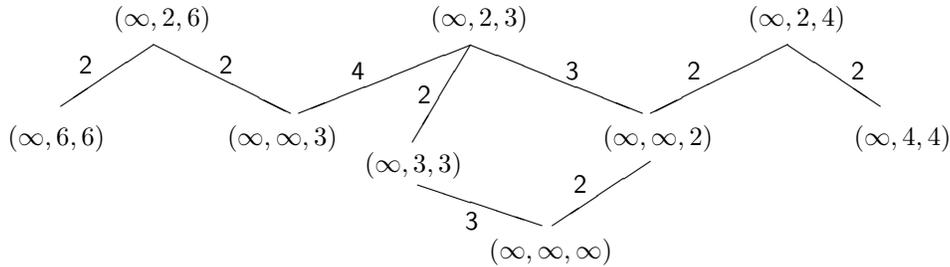
\begin{figure}[h!]
   \setlength{\unitlength}{1pt}
 \begin{picture}(150,250)(0,0)
               \put(90,240){$(\infty, 2, 6)$} \put(210,240){$(\infty, 2, 3)$} \put(330,240){$(\infty, 2, 4)$}
               \put(105,233){\line(-3,-2){35}} \put(105,233){\line(2,-1){52}} \put(225,233){\line(-5,-2){65}}
                \put(225,233){\line(5,-2){65}} \put(345,233){\line(3,-2){35}} \put(345,233){\line(-2,-1){52}}
                \put(225,233){\line(-3,-5){22}} \put(185,185){$(\infty, 3, 3)$}
                \put(50, 195){$(\infty, 6, 6)$} \put(133, 195){$(\infty, \infty, 3)$} \put(275, 195){$(\infty, \infty, 2)$}
                \put(370, 195){$(\infty, 4, 4)$} \put(205, 180){\line(3,-1){47}} \put(293, 189){\line(-3,-2){37}}
               \put(232,152){$(\infty, \infty, \infty)$}
               \put(77,222){2} \put(130,222){2} \put(180,220){4}  \put(205, 211){2} \put(261,220){3}  \put(307,220){2}  \put(369,220){2}
               \put(223,163){3}  \put(264,177){2}
 \end{picture}

 \vspace{-150pt}

\end{figure}

The $\mbox{}_2F_1$-functions for this diagram are said to be associated with elliptic curves,
modular forms, and elliptic integrals \cite{maillard}.  
They appear in many contexts.
To cover them, it suffices\footnote{Solutions related to entries of the diagram can be expressed in terms of
those three entries. However, the other entries can still be relevant if we want solutions of minimal size, see Section 5.3.3 (decompositions) in \cite{VJKThesis} for more.
Entry $(\infty, \infty, \infty)$ corresponds to writing solutions in terms of the elliptic integrals $K$ and $E$.}
to cover: $(k,\ell,m) = (\infty, 2, m)$ with $m \in \{3,4,6\}$.
Some parts of this paper focus on $(k, \ell, m) = (\infty, 2, 3)$, but our website \cite{dessindatabase}
covers $m=4$ and $m=6$ as well. 

%

\pagebreak 

\subsection{Project Outline}
\label{introduction}

\begin{defn} \label{def11} Let\, $f: \mathbb{P}^1\rightarrow \mathbb{P}^1$ be a rational function.
Let $k,\ell,m$ be positive integers or $\infty$, {see remark \ref{Infinity-Case}}.
The $(k,\ell,m)$-exceptional points of $f$ are:
\begin{itemize}
 \item roots of $f$ of order not divisible by $k$,
  \item roots of $1-f$ of order not divisible by $\ell$,
  \item roots of $\frac1f$ of order not divisible by $m$.
\end{itemize}
{We denote the set of $(k,\ell,m)$-exceptional points of $f$ as $\text{E}_{_{k\ell m}}(f)$, or simply as $\text{E}(f)$ when $k, \ell, m$ are fixed.}
\end{defn}

\begin{remark} \label{Infinity-Case}
  If $k = \infty$ then all roots of $f$ are $(k,\ell,m)$-exceptional points (and likewise for roots of $1-f$ resp. $\frac1f$ if $\ell$ resp. $m$ is $\infty$).
\end{remark}

\begin{remark} \label{Count-old}
  {Our definition~\ref{def11} differs slightly from that in \cite{hoeijheun} which defined exceptional points as \\
\mbox{} \hspace{10pt}  \{roots of $f$ of order $\neq k$\} $\cup$ \{roots of $1-f$ of order $\neq \ell$ \} $\cup$ \{roots of $\frac1f$ of order $\neq m$\} \\
 { which contains our $E(f)$.}} 
\end{remark}

\noindent The goals are to:
\begin{enumerate}
\item[(a)] construct a database \cite{dessindatabase} that, up to M\"obius-equivalence, contains all
rational functions with five $(\infty,2,3)$-exceptional points,
and likewise for $(\infty,2,4)$ and $(\infty,2,6)$. 
\item[(b)] prove completeness (Sections \ref{section3} -- \ref{belyi-2maps})
\item[(c)] give a fast method for the following problem: Given a field $k \subseteq \mathbb{C}$ and $\{q_1,\ldots,q_5\} \subset \mathbb{P}^1$,
 find every $F \in k(x)$ such that  $\text{E}(F) = \{q_1,\ldots,q_5\}$ (Section \ref{FindF})
\item[(d)] solve linear differential equations with 5 true singularities (Section \ref{motivation}).
\end{enumerate}

\noindent Our work continues the work (for $d=4$ exceptional points) in \cite{vidunas bound, hoeijheun}; the main novelties are:
 \begin{enumerate}
\item We prove completeness by giving an efficient new algorithm for finding dessins (algorithm 4.4). 
Such an algorithm was not needed for \cite{vidunas bound, hoeijheun};
the tables in \cite{vidunas bound} are small enough for manual enumeration, while
\cite{hoeijheun}
gave a method specific to $d=4$ that did not rely on dessins.
\item We compute almost-Belyi maps 
(\cite{vidunas bound,hoeijheun} only consider Belyi maps),
and braid orbits of almost-dessins
to prove completeness.
\item Our ${\rm Belyi}^{(1)}$ families turn out to be remarkably nice: they allow rational parametrizations that cover all ${\rm Belyi}^{(1)}$ maps by direct substitution, without any gaps or duplicates
(definitions 2.2 and 5.2).
\end{enumerate}


The first page on our website  \cite{dessindatabase}  gives the database for goal (a),
while another page (follow the link on the line ``Completeness'') gives examples and all algorithms needed for goals (b), 
(c) and (d).
For (c),
one has to select every ${f}$ in the database
whose $\text{E}(f)$ 
match $\{q_1,\ldots,q_5\}$ up to a M\"{o}bius transformation $x \mapsto \frac{ax+b}{cx+d}$.
We do this by computing
{\em five-point-invariants}
(functions of $\{q_1,\ldots,q_5\}$
whose values are invariant under M\"{o}bius transformations of the input\footnote{A
four-point invariant is given by the $j$-invariant of $y^2=(x-q_1)(x-q_2)(x-q_3)(x-q_4)$.}).

\begin{table}[H]
\begin{center}
    \begin{tabular}{|c|c|c|c|c|c|c|c|}
    \hline
   $(k, \ell, m)$ & \multicolumn{4}{|c|}{Belyi maps} & \multicolumn{2}{|c|}{${\rm Belyi}^{(1)}$ families} & {${\rm Belyi}^{(2)}$} \\ \cline{2-7}
    &{$\notin {\rm Belyi}^{(1)}$}& \multicolumn{2}{|c|}{$\in {\rm Belyi}^{(1)}$} & Total   &{$\notin {\rm Belyi}^{(2)}$}& {$\in {\rm Belyi}^{(2)}$} & families \\ \cline{3-4}
    && indirectly & directly &&&&\\ \hline
    $(\infty, 2, 3)$  & ${\it 411}$ &$9$& $266$& $686$		& ${\it 65}$ & $3$ & ${\it 2}$  \\ \hline
    $(\infty, 2, 4)$ & ${\it 121}$& $3$&$23$&$147$		&${\it 20}$&$0$&${\it 0}$\\ \hline
    $(\infty, 2, 6)$ & ${\it 54}$& $2$ &$5$&$61$			&${\it 12}$&$0$&${\it 0}$\\ \hline

    \end{tabular}
    \caption{Summary of the online table \cite{dessindatabase}.} \label{SummaryTable}
 \end{center}
 \end{table}

All maps in Table~\ref{SummaryTable} are listed on our website, and all have $|E(f)| = 5$. The three columns highlighted in italics
(Column ``$\not\in {\rm Belyi}^{(1)}$'', Column ``$\not\in {\rm Belyi}^{(2)}$'', and Column ``${\rm Belyi}^{(2)}$'')
contain precisely those entries of Table~\ref{SummaryTable} that would still have 5 exceptional points if we used the definition in Remark~\ref{Count-old}.
So although our application uses definition~\ref{def11},  both definitions are useful for the construction of the database.
The files on our website give algorithms and tables for both definitions, using the phrase ``count=5'' to refer to our definition~\ref{def11}, and the phrase ``Count=5'' for Remark~\ref{Count-old}.

All entries of the table are needed for goal (b), proving completeness. However, the only entries that are needed for goal (c) are the entries in italics plus Column ``indirectly'',
as will be explained in Section~\ref{BB1}.

Belyi maps and almost-Belyi maps have other applications as well. Indeed, some entries of our table have appeared elsewhere.
The almost-Belyi maps of degree $\leq$ 4 are constructed in \cite{AK}. In addition, 14 out of the 68 $\text{Belyi}^{(1)}$ maps corresponding to $(k, \ell, m) = (\infty, 2, 3)$ in Table~\ref{SummaryTable} are constructed in \cite{Kit1, Kit2, AK}.

Our main goal for this paper to prove completeness of our online table (summarized in Table~\ref{SummaryTable}) since that
will be useful for our application. 
The Belyi and $\text{Belyi}^{(1)}$ columns have so many entries that algorithms are needed for this proof.
Algorithms that are key to the proof will be described in this paper.
Particularly important is Algorithm~4.4 which is key to proving that our Belyi tables are complete.
It computes all ``dessins'' (conjugacy classes of 3-constellations) relevant for our project.
Although there is another implemented algorithm for the same task \cite{Sijsling} based on group theory,
it was not efficient enough {for higher degrees}, so we developed a novel algorithm instead.



\section{Rational functions with a prescribed branching pattern} \label{section2}
This section will cover goal~(a) from the introduction.
Section~\ref{branching_pattern} will enumerate the relevant branching patterns. Finding function(s) for a branching pattern, as in \cite{RamifiedCoverings,SijslingVoight},
is shown here by an example:

\begin{exmp}
\label{computeF}
Suppose we want to find a rational function $f: \mathbb{P}^1 \rightarrow \mathbb{P}^1$
with branching pattern  $(1,1,3,5)$, $(2,2,2,2,2)$, $(1,3,3,3)$ above $0,1,\infty$.
We abbreviate this as $(1^2, 3, 5), (2^5), (1, 3^3)$.
The degree is \linebreak $n =  1+1+3+5 = 2+2+2+2+2 = 1+3+3+3 = 10$.
The sum of $e_p-1$ (where $e_p$ denotes the branching index) for all $p$ above $\{0,1,\infty\}$
is $(1-1)+\cdots+(3-1) = 17$.  However, the Hurwitz formula for genus zero:
\begin{equation}  \sum_{p \in \mathbb{P}^1} (e_p - 1) = 2n - 2 \label{Hur} \end{equation}
gives 18, so there must be 1 more ramification point $p$, with $e_p = 2$, above some point $t \not\in \{0,1,\infty\}$.
So, $f$ is a {\em $\text{Belyi}^{(1)}$} map as it has a simple branch point outside $\{0,1,\infty\}$.

We aim to find all such $f$ up to M\"{o}bius-equivalence (Definition~\ref{me}). We use the three degrees of freedom
in M\"{o}bius-transformations to move the order-1 pole of $f$ to $x=1$,
and the roots of orders 3 and 5 to $x=0$ and $x=\infty$.
That brings $f$ in this form:
$$f:= \frac{(Ax^2+Bx+C)x^3}{(x-1)({x}^{3}+{ a_2}{x}^{2}+{ a_1
}x+{ a_0} )^{3}},  \,\,\,\, 1-f = {\frac {c \left( {x}^{5}+{ b_4}{x}^{4}+{ b_3}
{x}^{3}+{ b_2}{x}^{2}+{ b_1}x+{ b_0} \right) ^{2}}{\left( x-1 \right)  \left( {x}^{3}+{ a_2}{x}^{2}+{ a_1
}x+{ a_0} \right) ^{3}}}.$$
%
%
Equating $f$ with $1-(1-f)$ produces equations for the unknowns. 
Our implementation 
eliminates unknowns as long as it finds an equation that is linear in an unknown.
Three unknowns $b_2,b_3,b_4$ in two large equations remain. Factoring the resultant produces one equation in two unknowns,
i.e. an algebraic curve which turned out to have genus 0
(remarkably, the same happened for all $68+20+12$ cases in Tables~\ref{tableB1} and~\ref{table46}).
That means the solutions to this equation can be written as rational functions in a new variable $s$, which we can find with
Maple's {\tt parametrization}. After simplification we obtain
$b_3 = \frac13 s^4  - \frac83 s^3  + 18 s^2  - 96 s + 18$ and $b_4 = s^2-14$. Substitution followed by a gcd produces the value of $b_2$.  
Substituting into $f$, followed by two simple transformations ($x \mapsto 1-x$, and $s \mapsto s+3$) to reduce its size, produces:
$$f = \frac{64s^8(x-1)^3(9x^2-6s^2x-28sx+s^4-12s^3+36s^2)}{x (9x^3-6s^2x^2-36sx^2+s^4x-4s^3x+60s^2x+8s^4-32s^3)^3} \in \mathbb{Q}(s)(x).$$
This branches above $0,1,\infty$ plus one more point, denoted $t$.
To find it we  first compute the ramification point $p$ above $t$. This $p$ must be the only root of $f' = 0$
not in $f^{-1}(\{0,1,\infty\})$.
We find $p = \frac{1}{15} {s}^{2}-\frac23 s-\frac85$ and 
\begin{equation} \label{valuet}
	t = f(p) = {\frac {3125 \left( s-9 \right) ^{4}{s}^{8} \left( 3s
-2 \right) }{ 4\left( s-4 \right) ^{4} \left( {s}^{3}-9{s}^{2}+324s
-216 \right) ^{3} \left( s-1 \right) ^{2}}}.
\end{equation}
\end{exmp}

If $g: \mathbb{P}^1 \rightarrow \mathbb{P}^1$ has the same branching pattern 
above $0,1,\infty$,
one could ask if it is M\"obius-equivalent to $f$ for some value of $s$.
More generally, how to prove completeness for the entire database? Before we can answer that in Section~5 we first need a definition:

\begin{defn} \label{specialvalues} 
For a \emph{$\text{Belyi}^{(1)}$} map $f = f(s,x) \in \mathbb{C}(s)(x)$, let $\phi_f  \in \mathbb{C}(s)$ be the function that expresses
(as for example in equation~(\ref{valuet}))
the branch point $t \not\in \{0,1,\infty\}$ in terms of $s$. A point $s_0 \in \mathbb{P}^1$ is {\em degenerate} if $f(s_0, x)$ does not evaluate
to some $g \in \mathbb{C}(x)$ with the same $x$-degree as $f$.
It is called {\em generic} if $\phi_f(s_0) \not\in \{0,1,\infty\}$, and {\em special} if it is not degenerate nor generic. We define
$f$'s {\em family} as $\{f(s_0,x)\,|\,s_0$ not degenerate$\}$, and call it {\em gap-free} if no generic $s_0$ degenerates.
\end{defn}

\subsection{Enumerating branching patterns} 
\label{branching_pattern}

Let $B = (e_{1,1},\ldots,e_{1,n_1}), (e_{2,1},\ldots,e_{2,n_2}), (e_{3,1},\ldots,e_{3,n_3})$ be a {\em branching pattern of degree $n$},
which means that the $e_{i,j}$ are positive integers with $n = \sum_{j=1}^{n_i} e_{i,j}$ for each $i = 1,2,3$.

Since we only tabulate rational functions, we only consider {\em planar} (genus zero)
branching patterns. Then $S \leq 2n - 2$, where $S := \sum_{i=1}^3 \sum_{j=1}^{n_i} (e_{i,j} - 1)$
is the part of the Hurwitz formula~(\ref{Hur})
coming from points $p$ above $\{0,1,\infty\}$.
Let $\delta = 2n-2-S$. If $\delta=0$ then we call $B$ a {\em Belyi branching pattern},
if $\delta>0$ then we call $B$ a $\text{Belyi}^{(\delta)}$ branching pattern (Example~\ref{computeF} was planar and $\text{Belyi}^{(1)}$).
\begin{defn} \label{count}
Let
$B$ and $e_{i,j}$ as above
and $k,\ell,m$ be positive integers or $\infty$.
Let $(A_1,A_2,A_3) = (k, \ell, m)$. We define
{$\text{E}(B) := \{ (i,j)\, | \,\, A_i = \infty \,\, \text{or} \,\, A_i \,\, |\hspace{-8pt}\not \,\,\,\,  e_{i,j}\}$}. 
\end{defn}

\begin{remark} \label{Old-count}
{Definition \ref{count} corresponds to Definitions \ref{def11} in that if B is the branching pattern of a rational function $f$, then
  $|\, \text{E}(B)\,| = |\, \text{E}(f)\,|$.}
\end{remark}


If $B$ is planar $\text{Belyi}^{(\delta)}$ of degree $n$, then $\# e_{i,j} = \sum e_{i,j} - \sum (e_{i,j}-1) = 3n - S = n+2+\delta$.
The number of $e_{i,j}$ divisible by 
$A_i$ is at most $n/k + n/\ell + n/m$. So
if $(k, \ell, m) = (\infty, 2, 3)$ and $d = |\,\text{E}(B)\,|$, then $d \geq n+2+\delta - (n/\infty + n/2 + n/3)$ and hence
\begin{equation}
	n \leq 6(d-2-\delta). \label{boundn}
\end{equation}

Our website \cite{dessindatabase} has a routine (similar to \cite[Section~3]{hoeijheun}) to enumerate the necessary branching patterns.



%

\section{Riemann existence theorem and (almost)-Belyi maps}
\label{section3}
\begin{defn} \label{me} Two rational functions $f,g: \mathbb{P}^1 \rightarrow \mathbb{P}^1$ 
are called\, {\em M\"{o}bius-equivalent} if 
$f = g \circ m$ for some $m \in {\rm Aut}(\mathbb{P}^1)$ (=
the group of M\"{o}bius transformations\, $\{\frac{ax+b}{cx+d}\, |\, ad-bc \neq 0\}$). 
\end{defn}


\begin{defn}  $\cite{constellation}$
  A list $[g_1,\ldots,g_k]$ of permutations in $S_n$ is called a {\em $k$-constellation} if
    the group $\langle g_1,\ldots,g_k\rangle$ acts transitively on $\{1\ldots n\}$ and
    $g_1 \cdots g_k = 1$.
   Here 
   $n$ is the {\em degree}, and $\langle g_1,\ldots,g_k\rangle$ is the \emph{monodromy group} of $[g_1,\ldots,g_k]$.
\end{defn}
\begin{defn}
  $[g_1,\ldots,g_k]$ and $[h_1,\ldots,h_k]$ are
{\em conjugated} if $\exists_{\tau \in S_n} \forall_i \,\, h_i = \tau^{-1} g_i \tau$.
We denote the conjugacy class of $[g_1, \ldots, g_k]$ as $[g_1, \ldots, g_k]_\sim$.
  \end{defn}
  \begin{theorem}{\textbf{Riemann Existence Theorem}} (formulation from~\cite{RETNotes}, for more see \cite{RETBook,RETPaper1,RETPaper2}).
  Let $p_1,\ldots, p_k$ be distinct points of $\mathbb{P}^1$. For any transitive 
  representation $\rho: \pi_1(\mathbb{P}^1\setminus \{p_1,\ldots,p_k\}) \rightarrow S_n$ there is a connected Riemann surface $X$ and a proper holomorphic map $f: X \rightarrow \mathbb{P}^1$
  of degree $n$ which realizes $\rho$ as its monodromy homomorphism. Moreover $X$ and $f$ are unique up to equivalence.
 \end{theorem}

 \begin{remark} \label{Fieldk}
 If $p_1, \ldots, p_k \in k \cup \{\infty\}$ for a subfield $k \subseteq \mathbb{C}$ then $f$ can be defined over some algebraic extension of $k$ (Cor. 7.10 in \cite{RETBook}).
In Example~\ref{computeF}, $p_1,\ldots,p_4 \in k \bigcup \{\infty\}$ where $k := \mathbb{Q}(t)$,
while $f$ is defined over $\mathbb{Q}(s)$, 
an algebraic extension of $k$ (equation~(\ref{valuet}) shows $t \in \mathbb{Q}(s)$). 
 \end{remark}

 If the branched set $\{p_1,\ldots,p_k\}$ is $\{0, 1, \infty\}$ then the pair $X,f$ is called a {\em Belyi map}. 
 The representation $\rho$ is given by a $k$-constellation $[g_1,\ldots,g_k]$.
 We use this for Belyi and $\text{Belyi}^{(1)}$ maps in Sections~4 and~5, but not for $\text{Belyi}^{(2)}$ maps where we only have 2 cases (Section~6).
 \[ p_1,\ldots,p_k
  = \left\{ \begin{array}{l}
 	0,1,\infty 	\hspace{28pt} (k=3, \textrm{ Belyi\ case),} \\
 	0,1,t,\infty		{\rm \ \ \ \ \ \ } (k=4\,\ \textrm{where} \ g_t \ \textrm{is a 2-cycle, $\text{Belyi}^{(1)}$ case).} \end{array} \right. \]
%
We only use {\em planar} $k$-constellations, i.e. $X=\mathbb{P}^1$. 
Then $[g_1,\ldots,g_k]_\sim$ determines $f$ up to Aut($\mathbb{P}^1$):
\begin{equation}
	[g_1,\ldots,g_k]_\sim \, \Longleftrightarrow \, f \text{ up to M\"{o}bius-equivalence.}
	\label{fcor}
\end{equation}
In the Belyi case $[g_1,\ldots,g_k]_\sim$ corresponds to a {\em dessin d'enfant} as well.


\subsection{Dessins d'enfants} \label{dessins}
%

 \begin{defn}
 A {\em dessin d'enfant} \cite{LS,constellation} is a connected and oriented graph with black and white vertices,
where any edge joins a black and a white vertex. The {\em degree} is the number of edges.
 \end{defn}
 \begin{figure}[h!]
 \vspace{-30pt}
\setlength{\unitlength}{1pt}
 \begin{picture}(150,110)(-10,40)

 \put(45,73){$\bullet$} \put(48,75){\line(1,0){18}} \put(65,73){$\circ$}\put(69,75){\line(1,0){18}} \put(85,73){$\bullet$}
 \curve(48,75,32,88,17,77) \put(14,73){$\circ$} \curve(47,74,32,61,17,74)
 \curve(87,74, 89,82, 97,92) \put(96,91){$\circ$}\curve(100,94, 110,100, 120,103)
 \curve(87,74,89,66, 97,56) \put(96,53){$\circ$} \curve(100,54,110,48,120,45)
  \put(119,101){$\bullet$}
 \curve(121,103, 126,113, 139,120)\put(138,118){$\circ$} \curve(142,120, 155, 113, 160,104)
 \curve(121,103, 126,93, 139, 86) \put(138,84){$\circ$} \curve(142, 86, 155,93, 160, 104) \put(158,103){$\bullet$}
 \put(119,42){$\bullet$}
 \curve(121,44,126,54,139,61) \put(138,59){$\circ$} \curve(142,61,155,54,160,45)
  \curve(121,44, 126,34, 139, 27) \put(138,25){$\circ$} \curve(142, 27, 155,34, 160, 45) \put(158,44){$\bullet$}
 \curve(160,106, 174, 98, 186,77) \put(184,73){$\circ$}
 \curve(160,46 , 174,54, 186,74)
  \put(31,52){1} \put(30,79){3} \put(56,77){2} \put(74,66){4} \put(92,63){5} \put(84,84){6} \put(103,40){7} \put(106,90){10} \put(128,92){11}
  \put(118,117){12} \put(143,107){18} \put(155, 85){16} \put(177,96){17} \put(169, 62){14} \put(157,30){13} \put(130,32){8} \put(122,57){9} \put(144,48){15}
  \put(260,75){$\circ$}\put(264,77){\line(1,0){20}} \put(283,75){$\bullet$}
  \curve(285,77, 292,92, 305,101) \put(304,99){$\circ$} \curve(308,102, 320, 102, 331, 97 )
  \curve(285,77, 292,62, 305,53) \put(304,50){$\circ$} \curve(308,52, 320, 52, 331, 57) \put(330,56){$\bullet$}
  \put(330,94){$\bullet$} \curve(331,57, 338, 65, 342, 77) \put(339,76){$\circ$} \curve(333, 96, 339,89, 342, 80)\curve(332,59, 340, 50, 348, 41)
  \curve(332,96, 340, 104, 348, 112)  \put(347,111){$\circ$} \put(347,38){$\circ$}
  \put(273,68){1}\put(294,64){2}\put(286,92){3}\put(317,92){4} \put(341,88){5}\put(334,104){6}\put(317,43){7}\put(342,49){8} \put(332,68){9}
\end{picture}
\vspace{9pt}
\caption{two planar dessins d'enfants (the labels are not part of the dessins d'enfants)}
\label{planardessins}
\end{figure}
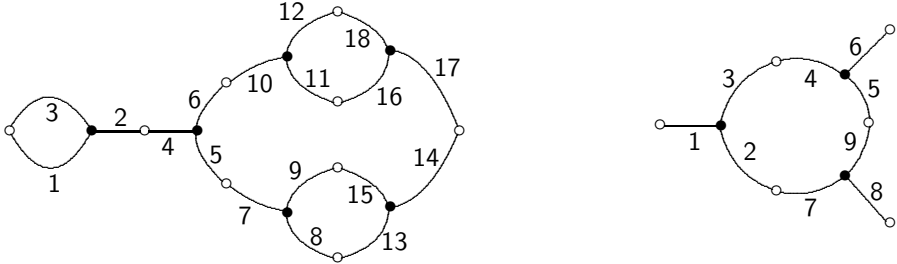

\noindent There is a one-to-one correspondence \cite{beukers}  between:
\begin{enumerate}
\item 3-constellations up to conjugation,
\item Belyi maps up to equivalence,
\item dessins d'enfants up to equivalence.
\end{enumerate}
$1 \mapsto 2$: The Riemann existence theorem. \\[5pt]
$2 \mapsto 3$: 
The dessin d'enfant of a Belyi map $f$ is the graph $f^{-1}([0,1])$ where:
 $f^{-1}(\{0\}) = \{$black vertices$\}$,
 $f^{-1}(\{1\}) = \{$white vertices$\}$, and
 $f^{-1}(\,(0,1)\,) = \{$edges$\}$.
 Faces correspond 1-1 
 to $f^{-1}(\{\infty\})$. \\
Example: $f_1 = 4 (x^6-4x^5+5x^2+4x+4)^3/(27 (x-4) (5x^2+4x+4)^2 x^5)$ and $f_2 = 4(x^3-1)^3/(27x^3)$.
Plotting $f_1^{-1}([0,1])$ and $f_2^{-1}([0,1])$ we find Figure~\ref{planardessins} up to homeomorphism, without the labels \cite{dessindatabase}.
\\[5pt]
\noindent $2 \mapsto 1$: The monodromy  command in Maple
computes a $k$-constellation for any algebraic function. Applying this (see Section 3.1 in \cite{dessindatabase}) to $f_1$ produces: \\
\mbox{} \,\,\,  $g_0 = (1\ 5\ 3)(2\ 4\ 6)(7\ 9\ 11)(8\ 12\ 10)(13\ 14\ 15)(16\ 18\ 17)$ \\
\mbox{} \,\,\,  $g_1 = (1\ 4)(2\ 8)(3\ 7)(5\ 9)(6\ 12)(10\ 13)(11\ 14)(15\ 16)(17\ 18)$ \\
\mbox{}\hspace{8pt}$g_{\infty} = (1\ 2\ 10\ 15\ 17\ 16\ 14\ 9)(3\ 11\ 13\ 12\ 4)(5\ 7)(6\ 8)$. \\[5pt]

\pagebreak 

$3 \mapsto 1$:
After adding a label to each edge (the labels are not part of the dessin d'enfant itself)
we can read the 3-constellation $[g_0,g_1,g_{\infty}]$ from the ``labelled dessin'' as follows.
Reading labels counter-clockwise around each black vertex
produces each cycle of $g_0$ (some of which may be 1-cycles, the valence of a vertex is the length of the corresponding cycle).
Likewise, each white vertex corresponds to a cycle of $g_1$.
There are two ways to find $g_{\infty}$, one could compute it as $(g_0 g_1)^{-1}$,
but one can also read $g_{\infty}$ directly from the ``labelled dessin'';
each cycle of $g_{\infty}$  is found by following the labels inside each face. From the first ``labelled dessin'' we read: \\
\mbox{} \,\,\,$h_0 \, = {(1\ 2\ 3)(4\ 5\ 6)(7\ 8\ 9)(10\ 11\ 12)(13\ 14\ 15)(16\ 17\ 18)}$ \\
\mbox{} \,\,\,$h_1 \, = {(1\ 3)(2\ 4)(5\ 7)(6\ 10)(8\ 13)(9\ 15)(11\ 16)(12\ 18)(14\ 17)}$ \\
\mbox{} \,\,$h_{\infty} = {(3)(1\ 2\ 6\ 12\ 17\ 13\ 7\ 4)(5\ 9\ 14\ 16\ 10)(8\ 15)(11\ 18)}$. \\
Algorithm~4.3 in Section~\ref{ComputingDessins} can verify that the 3-constellations $[g_0,g_1,g_{\infty}]$ and $[h_0,h_1,h_{\infty}]$ are conjugated.
Several algorithms in Section~\ref{section4} use permutations in {\em expanded form}, which means the 1-cycles are written as well. For example, the 3-constellation of the second ``labelled dessin'' is: \\
\mbox{} \,\,\, $[(1\ 2 \ 3)(4 \ 5 \ 6)(7 \ 8\ 9), \ (1)(6)(8)(2 \ 7)(3 \ 4)(5 \ 9), \ (1\ 7\ 8\ 5 \ 6 \ 3)(2\ 4 \ 9)]$.



\begin{remark}
  {If $D$ is the dessin of a Belyi map $f$ with branching pattern $B$, then
  $|\, E(f)\,|$,  $|\, E(B)\,|$, $|\, E (D)\,|$ denote the number of exceptional: points of $f$ resp. branchings in $B$ resp. cycle-lengths in $D$.
  Since these numbers are equal, we will also use the shorter notation $|E|$ if $f$, $B$, or $D$ is clear from the context.}
\end{remark}

\section{Belyi maps} \label{section4}
The goal in this section is Algorithm 4.4 which can compute
all dessins with $|\,\text{E}\,| = 5$, which is the key step to proving that
all Belyi maps $f$ with $|\,\text{E}\,| = 5$ appear in our table.

\subsection{Computing 3-constellations} \label{compute3-constln}

\begin{defn}
Let $g \in S_n$.  Then $g'$ denotes an element of $S_{n-1}$ defined as follows: for $i \in \{1\ldots n-1\}$ define $g'(i)$ as $g(i)$ if $g(i) \neq n$,
and $g(g(i))$ if $g(i) = n$. When $g$ is written in disjoint cycle notation, one obtains $g' \in S_{n-1}$ by simply erasing $n$.
\end{defn}

\begin{figure}[H]
\setlength{\unitlength}{1pt}
 \begin{picture}(400,250)(0,0)
 \put(20,240){$T_1$:} \put(180,240){$\bullet$} \put(185,242){\line(1,0){30}} \put(214,240){$\circ$} \put(200,234){\small 1} \put(228, 244){$g_0 = (1)$}
                        \put(228,231){$g_1 = (1)$}  \put(200,222){\vector(-2,-1){95}} \put(200,222){\vector(0,-1){43}} \put(200,222){\vector(3,-1){100}}
 \put(20, 165){$T_2$:} \put(65,165){$\circ$} \put(70,167){\line(1,0){20}} \put(89,165){$\bullet$} \put(93,167){\line(1,0){20}}
                       \put(113,165){$\circ$} \put(78, 170){\small 2} \put(103,159){\small 1} \put(70, 146){$g_0= (1\,2)$} \put(70,133){$g_1 = (1)(2)$}
                         \put(175,165){$\bullet$} \put(180,167){\line(1,0){20}} \put(200,165){$\circ$} \put(204,167){\line(1,0){20}} \put(223,165){$\bullet$}
                         \put(188, 159){\small 1} \put(212,170){\small 2} \put(177, 146){$g_0= (1)\,(2)$} \put(177,133){$g_1 = (1 \, 2)$}
                        \curve(290,170,310,185,330,170) \put(287,166){$\bullet$} \put(328,166){$\circ$} \curve(290,167, 311,152, 330,167)
                         \put(310, 144){\small 1} \put(310,175){\small 2}  \put(338, 169){$g_0 = (1\, 2)$} \put(338,153){$g_1 = (1 \, 2)$}
                         \put(90,128){\vector(0,-1){42}} \put(90,128){\vector(1,-1){48}} \put(200,126){\vector(0,-1){42}} \put(200,126){\vector(-1,-1){45}}
                         \put(200,126){\vector(1,-1){45}}   \put(310,136){\vector(0,-1){27}} \put(310,136){\vector(-3,-4){42}}
 \put(20, 70){$T_3$:} \put(65,70){$\bullet$} \put(69,72){\line(1,0){11}} \put(79,70){$\circ$} \put(83,72){\line(1,0){11}}
                          \put(93,70){$\bullet$} \put(97,72){\line(1,0){11}} \put(108,70){$\circ$} \put(101,64){\small 1} \put(87,74){\small 2} \put(73,64){\small 3}
                          \put(60,48){$g_0=(1\,2)(3)$} \put(60,35){$g_1= (1)(2\,3)$} \put(125,70){$\cdots \cdots \cdots$}
                           \put(179,70){$\circ$} \put(183,72){\line(1,0){11}} \put(193,70){$\bullet$} \put(197,72){\line(1,0){11}}
                          \put(207,70){$\circ$} \put(211,72){\line(1,0){11}} \put(221,70){$\bullet$} \put(215,74){\small 2} \put(201,64){\small 1} \put(187,74){\small 3}
                          \put(174,48){$g_0 = (1)(2\,3)$} \put(174,35){$g_1 = (1\,2)(3)$} \put(242,70){$\cdots \cdots \cdots$}
                           \curve(295,75,325,85,345,75) \put(292,71){$\bullet$} \put(343,71){$\circ$} \curve(295,72, 321,47, 345,72)
                         \put(305, 43){\small 1} \put(308,74){\small 2} \put(306,93){\small 3} \put(353, 75){$g_0=(1\, 2\,3)$} \put(353, 62){$g_1=(1\,2\,3)$}
                          \curve(295,75, 307,103, 320,87)
                         \curve(320,83, 327,65, 343,71)
\end{picture}
\vspace{-30pt}
\caption{Computing repeat-transitive 3-constellations (definition 4.2)}
\label{computing3constellations}
\end{figure}
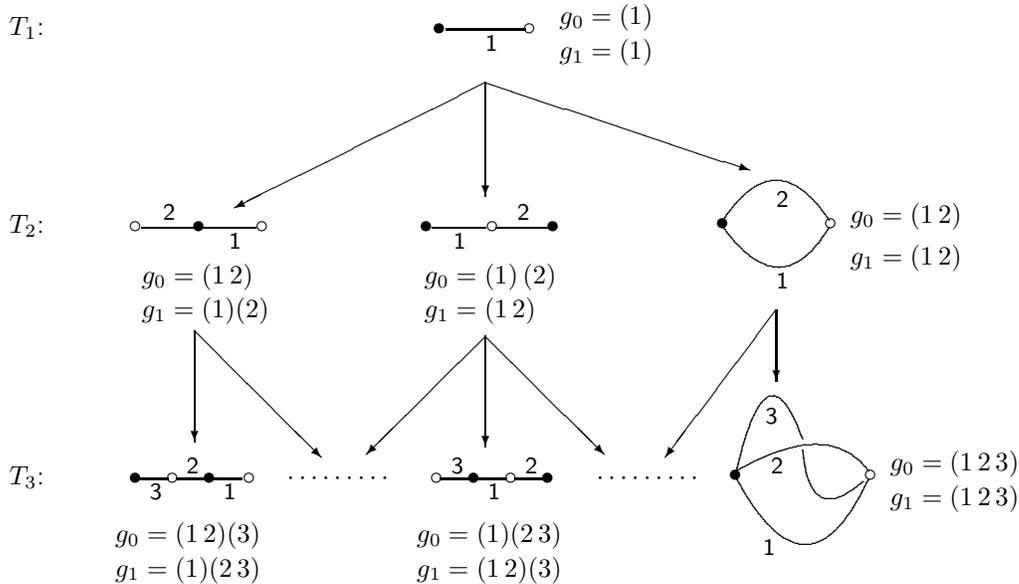

\begin{defn}
 A pair $(g,h)$ where $g,h \in S_n$ is called repeat-transitive when $n=1$, or, when $<g,h>$ is a transitive subgroup of $S_n$ and $(g',h')$ is repeat-transitive in $S_{n-1}$.
\end{defn}

Figure~\ref{computing3constellations} shows (for $N=3$) how one can compute all repeat-transitive 3-constellations of degrees $1,2,\ldots,N$ recursively.
Start with the 3-constellation of degree 1. Then, given the set $T_{n-1}$ of \linebreak repeat-transitive
3-constellations of degree $n-1$,
insert one more edge in all possible ways to obtain all repeat-transitive 3-constellations of degree $n$.
Algorithm~4.1 below shows how to implement this. 
It represents 3-constellations as $[g_0,g_1]$ (since $g_{\infty}$ can be recovered from $g_0 g_1 g_{\infty} = 1$)
with $g_0,g_1 \in S_{n-1}$ written in {\em expanded form}, i.e., including 1-cycles.
Inserting an edge means doing the following to both $g_0$ and $g_1$:
  \begin{enumerate}
  \item[(i)] insert the new number $n$ into an existing cycle, or
  \item[(ii)] add a new 1-cycle $(n)$.
  \end{enumerate}
One can not choose (ii) for both $g_0$ and $g_1$, because the resulting pair would not be transitive 
(the graph would not be connected).
This leaves $n^2-1$ ways in Step~2 of Algorithm~4.1 to add an edge to $[g_0,g_1]$. Indeed, $|T_2| = 3 = (2^2-1) \cdot |T_1|$ in Figure~\ref{computing3constellations}.
In Step~2,
the call {\tt Insert}$(g_0, i,n)$ (with $g_0 \in S_{n-1}$ and $i \in \{1,\ldots,n\}$) inserts edge $\#n$ at the $i^{th}$ position, as shown here:
\begin{exmp}
Let  $g_0 = (1\, 2) (4\, 5) (6\, 8) \in S_8$. The program call {\tt Insert}$(g_0,6,9)$ computes: \\ 
{\rm \textbf{Step 1:}} Write $g_0$ in expanded form (including $1$-cycles) so that all edges 1--8 appear. Placeholders (asterisks) indicate all $9$ possible positions in $g_0$ where $9$ can be inserted:
$$g_0 = (3\, \ast) (7\, \ast) (1 \ast 2\,\ast) (4\ast 5\,\ast) (6\ast 8\,\ast) (\ast)$$
{\rm \textbf{Step 2:}} Insert 9 in the $6^{th}$ placeholder:
$${\tt Insert}(g_0, 6, 9):= (3) (7) (1\, 2) (4\, 5\, 9) (6\, 8) \ \ \ {\rm (written \ in \ expanded \ form).}$$
\end{exmp}

\begin{framed}
\noindent \textbf{Algorithm~4.1: Compute all repeat-transitive 3-constellations of degree $\leq N$}.  \\
\textbf{Step 1:} $T_1 := \{[ (1), (1) ] \}$  \\
\textbf{Step 2:} For $n$ from 2 to $N$ do:\\
\mbox{} \hspace{37pt}$T_n := \{[{\tt Insert}(g_0, i, n), {\tt Insert}(g_1, j, n)]\,\,|\,\,i,j \in \{1\ldots n\},\, (i,j)\neq(n,n),\, [g_0, g_1] \in T_{n-1}\}.$
\end{framed}

 From $|T_1| = 1$ and $|T_n| = (n^2-1)\cdot |T_{n-1}|$ one finds
 $$|T_n| =  {(n-1)! (n+1)!}\hspace{1pt}/\hspace{1pt}2 = 1,\, 3,\, 24,\, 360,\, 8640,\, 302400,\, 14515200,\, 914457600,\, \ldots $$
 \begin{remark} \label{notTR}
There are twenty-six 3-constellations of degree 3. Two of them are not repeat-transitive: \\[-15pt]
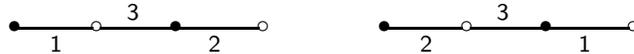
\begin{figure}[h!]
 \setlength{\unitlength}{1pt}
 \begin{picture}(50,20)
   \put(70,10){$\bullet$} \put(72,12){\line(1,0){30}} \put(101,10){$\circ$} \put(86,3){\small 1}
   \put(105,12){\line(1,0){30}} \put(131,10){$\bullet$} \put(135,12){\line(1,0){30}} \put(164,10){$\circ$} \put(115,15){\small 3}  \put(146,3){\small 2}

   \put(210,10){$\bullet$} \put(212,12){\line(1,0){30}} \put(241,10){$\circ$} \put(226,3){\small 2}
   \put(245,12){\line(1,0){30}} \put(271,10){$\bullet$} \put(275,12){\line(1,0){30}} \put(304,10){$\circ$} \put(255,15){\small 3}  \put(286,3){\small 1}

  \end{picture}

  \vspace{-10pt}
\end{figure}

\noindent So the set $T_3$ computed by Algorithm~4.1 has twenty-four $3$-constellations.
\end{remark}

The construction in Figure \ref{computing3constellations} is complete up to re-labeling (e.g., compare the two 3-constellations from Remark~\ref{notTR}
with $T_3$ in Figure \ref{computing3constellations}). So Algorithm~4.1 does find all 3-constellations up to conjugation. 
Dessins with {$|\,\text{E}\,| = 5$} have degrees $\leq 18$, see inequality~(\ref{boundn}). To find them, we need to implement several improvements because $T_{18}$ is much too large for the computer.

 \subsection{Discarding unnecessary 3-constellations}
 \label{DiscardNonPlanarDessins}

Let $\#g_0$ denote the number of cycles in $g_0$, including 1-cycles.
For a 3-constellation $[g_0,g_1,g_{\infty}]$ of degree~$n$, the genus $g$ of the corresponding dessin d'enfant is given by Euler's formula:
$$ 2-2g = \#{\rm vertices} - \#{\rm edges} + \#{\rm faces} = \#g_0 + \#g_1 - n + \#g_{\infty}.$$
Our aim is rational functions, which correspond to planar (i.e. $g=0$) 3-constellations.
Adding edges to a non-planar dessin d'enfant can not 
make it planar, so
we may discard non-planar 3-constellations in Algorithm~4.1 as soon as they occur.
This reduces the growth of $T_n$ but more
improvements are needed since it still grows much too fast.

The goal 
is {$|\,\text{E}\,| = 5$}, however, we can not simply discard 3-constellations with {$|\,\text{E}\,| > 5$} as soon as they occur,
because adding an edge can lower the value of 
$|\,\text{E}\,|$. We solve this problem with {\em weighted} counts.

\begin{defn} \label{wcount} Notations as in Section~\ref{branching_pattern}, Definition~\ref{count}. \newline
If $A_i | e_{i,j}$ then $s_{i,j} := 0$, if $e_{i,j} \equiv -1$ mod $A_i$ then $s_{i,j} := \frac12$, otherwise $s_{i,j} := 1$ (if $A_i = \infty$ then $s_{i,j} = 1$).
The weighted-count of $B$ is the sum of the $s_{i,j}$.
\end{defn}
\noindent If we replace $\frac12$ in Definition~\ref{wcount} by 1, then we get $|\,\text{E}\,|$ 
from Definition~\ref{count}. Thus,
\begin{equation} \mbox{\rm $|\,\text{E}\,| \, \geq \mbox{\rm weighted-count}$} 
\label{ineq} \end{equation}
\begin{proposition}  \label{p} Adding an edge does not decrease the weighted-count if the dessin d'enfant stays planar and $A_3 = \infty$.
\end{proposition}
\noindent {\bf Proof:}
Let $S_i$ be the sum of the $s_{i,j}$, and let $\tilde{S}_i$ be the sum after adding one edge. Then $\tilde{S}_1 = S_1 + 1$ if the number of black vertices increased, otherwise, $\tilde{S}_1 \geq S_1 - \frac12$.
Likewise,  $\tilde{S}_2 = S_2 + 1$ if the number of white vertices increased, otherwise, $\tilde{S}_2 \geq S_2 - \frac12$.
The number of faces is $S_3$ if $A_3=\infty$, and
does not decrease when adding an edge ($\tilde{S}_3 \geq S_3$) if the result remained planar.
Now $\tilde{S}_1 + \tilde{S}_2 + \tilde{S}_3 \geq S_1 + S_2 + S_3$ because if no black or white vertices were added, then adding an edge increases the number of faces. \eindebewijs

We now switch from $(\infty,2,3)$ to $(3,2,\infty)$, an easily reversible transformation (for Belyi maps it means $f \mapsto 1/f$).
Then we may discard 3-constellations with 
weighted-count $>5$ in Algorithm~4.1 as soon as they occur; adding edges 
can not lead to $|\,\text{E}\,| \leq 5$ by Proposition~\ref{p} and inequality~(\ref{ineq}).
This drastically reduces the growth of $T_n$, but a problem still remains, which we handle next.

 \subsection{Finding a unique representative of a conjugacy class}
 \label{ComputingDessins}
The 3-constellations $[(1\ 2)(3), (1)(2\ 3)]$ and $[(1)(2\ 3), (1\ 2)(3)]$ in Figure~\ref{computing3constellations} are conjugated. 
We should remove all but one constellation in each conjugacy class,
not only because this gives an another drastic reduction in the growth of $T_n$, but also because we need 3-constellations {\em up to conjugacy} for
the correspondence from Section~\ref{section3}. 
For $\tau \in S_n$, denote $g^{\tau} := \tau^{-1} g \tau$. 

\begin{framed}
\noindent \textbf{Algorithm 4.2: Sort With Base point (SWB)} \\
\textbf{Input:} Transitive $g_0\ldots g_s$ in $S_n$ and a base point $b\in \{1, \ldots, n\}$. \\
\textbf{Output:} $[g_0^\tau \ldots g_s^\tau]$ for some $\tau \in S_n$ with the \textbf{property}: $[g_0 \ldots g_s]$ is conjugated to $[h_0 \ldots h_s]$ if\\
\mbox{}\hspace{42pt} and only if $\{{\tt SWB}(g_0 \ldots g_s, b)\,|\,1\leq b \leq n\}\, = \, \{{\tt SWB}(h_0 \ldots h_s, b)\,|\,1 \leq b \leq n\}$. \\
\textbf{Step 1:} $\pi_1 := b$. \\
\textbf{Step 2:} For $k$ from 1 to $n-1$ let $\pi_{k+1} := g_i(\pi_l)$ where $(i,l)$ is the first pair in $\{0 \ldots s\} \times \{1 \ldots k\}$  \\
\mbox{}\hspace{36.5pt} with $g_i(\pi_l) \not\in \{\pi_1,\ldots,\pi_{k}\}$. \\
\textbf{Step 3:} Let $\tau \in S_n$ with $\tau(i)=\pi_i$ and return $[g_0^\tau \ldots g_s^\tau]$.
\end{framed}
Verifying that ${\tt SWB}(g_0 \ldots g_s,b) = {\tt SWB}(g_0^\tau \ldots g_s^\tau, \tau^{-1}(b))$ for $\tau \in S_n$ proves the claimed property.
\begin{framed}
\noindent \textbf{Algorithm 4.3: UniqueRepresentative} \\
\textbf{Input:} Transitive $g_0\ldots g_s$ in $S_n$. \\
\textbf{Output:} A unique representative in the $S_n$-conjugacy class of $[g_0 \ldots g_s]$. \\
\textbf{Step 1:} $S:=\,\{{\tt SWB}(g_0 \ldots g_s,b)\,|\, 1 \leq b \leq n\}$. \\
\textbf{Step 2:} Return the first (we use a lexicographic ordering) element of $S$.
\end{framed}


\subsection{Computing dessins to prove that the Belyi table is complete}
\label{4.4}
From here on, the phrase ``dessin'' is short for ``conjugacy class of $3$-constellations'' represented by the output of Algorithm~4.3.
The bound $6(d-2)$ comes from Equation (3).
\begin{framed}
\noindent \textbf{Algorithm 4.4: Compute all planar dessins with $|\,E\,| = d$} \\
\textbf{Step 1:} $T_1 := \{[ (1), (1) ] \}$  \\
\textbf{Step 2:} For $n$ from 2 to $6(d-2)$ do:\\
\mbox{} \hspace{37pt}$T_n := \{[{\tt Insert}(g_0, i, n), {\tt Insert}(g_1, j, n)]\,\,|\,\,i,j \in \{1\ldots n\},\, (i,j)\neq(n,n),\, [g_0, g_1] \in T_{n-1}\}$ \\
\mbox{} \hspace{37pt}$T_n := \{[g_0,g_1] \in T_n\,|\,[g_0,g_1,(g_0 g_1)^{-1}]$ is planar and has weighted-count $\leq d\}$\\
\mbox{} \hspace{37pt}$T_n := \{$UniqueRepresentative$(g_0,g_1)\,\, | \,\, [g_0, g_1] \in T_n\}$. \\
\textbf{Step 3:} Return $\{[g_0,g_1]\,\,|\,\,n \leq 6(d-2) {\rm \ and \ } [g_0,g_1] \in T_n\,\text{with}\, |\,\text{E}\,| = d\}.$
\end{framed}

Algorithm 4.4 produces the following dessins for $d \in \{4, 5\}$: 
\begin{table}[H]
\begin{tabular}{| l | l |}
  \hline  \cline{1-2}
 $d$ & Number of planar dessins with $(3,2,\infty)$-count $d$. \\
  \hline  \cline{1-2}
 4  & 0, 1, 3,\, { 5},\, 3,\, { 10},\, 4, 6, 4, 4, 0, 6 \\
 \hline \cline{1-2}
 5  & 0, 0, 2,\, { 10}, { 18}, { 40}, { 50}, { 71}, { 76}, { 103},\, 36,\, { 108},\, 40, 42, 32, 32, 0, 26 \\
\hline \cline{1-2}
\end{tabular}
\end{table} \vspace{-5pt}

Another way to generate maps, using parenthesis systems, was given in \cite{Walsh}.
{Although we are mainly interested in $d = 5$, we ran Algorithm 4.4 for $d \leq 7$,
for both definitions {(see Section~\ref{introduction})}.
The entries of degree $n=6(d-2)$ form sequence $2,6,26,191,1904,\ldots$ (\url{www.oeis.org/A112948}).
Beukers and Montanus \cite{beukers} computed the dessins and Belyi maps for $d=6$, $n=24$
with a combination of machine and hand computation, but at the time they missed one of the 191 dessins.
To avoid the likelihood of a gap in a large table, it is important to verify it with machine-only computation.} \\

\noindent {\bf Proving Completeness for the Belyi table:}
To prove that our website \cite{dessindatabase} lists all rational Belyi maps
with $|\,\text{E}\,| = 5$, it is not enough
to check that its number of functions of degree $n$ matches row $d=5$ in the above table.
That would leave open the possibility 
of M\"obius-equivalent (Definition~\ref{me}) duplicates while missing other ones.
So we implemented a more rigorous check \cite{dessindatabase}.
It computes the dessin for each $f$ in our table by applying $2 \mapsto 1$ from Section~\ref{dessins}, followed by Algorithm~4.3.
Completeness is then proven by comparing these dessins with the independently-computed set of dessins from Algorithm~4.4.\\[5pt]
{Invariants offer a much faster way to prove completeness of our Belyi table, without the time-consuming computation $2 \mapsto 1$ from Section~\ref{dessins}.
Each time a pair $f_1 \neq f_2$ in our Belyi table had the same branching pattern, it turned out that their five point invariants are not equal, i.e; $I_5(f_1) \neq I_5(f_2)$. More details about these invariants are given in Section \ref{FindF}.
This proves that the table has no M\"obius-equivalent duplicates. To prove completeness it now suffices to compare (for each branching pattern)
the number of Belyi maps in the table with the number of dessins from Algorithm~4.4.} \eindebewijs

{Our $(\infty,2,3)$-Belyi table~\cite{dessindatabase} has 255+9+99 = 363 entries
representing 411+9+266 = 686 functions $F^B_1,\ldots,F^B_{686} \in \mathbb{C}(x)$
(If $f \in \mathbb{Q}(\alpha)(x)$ with $[\mathbb{Q}(\alpha):\mathbb{Q}]=d$, then it
represents $d$ elements of $\mathbb{C}(x)$, one for each of the $d$ complex roots of the minimal polynomial of $\alpha$.)
Their dessins are precisely the 686 dessins produced by Algorithm 4.4
(there are no obstruction issues as in \cite[Section~6]{hoeijheun}).}


\section{$\text{Belyi}^{(1)}$ maps} \label{belyi-1maps}

We consider planar 4-constellations $[g_0,g_1,g_t,g_{\infty}]$ where $g_t$ is a 2-cycle.
The phrase ``almost-dessin'' in this section refers to: conjugacy class of such a 4-constellation, represented by the output of Algorithm~4.3.
Almost-dessins corresponds to $\text{Belyi}^{(1)}$ maps up to M\"obius-equivalence, see (\ref{fcor}) in Section~3.

\subsection{Finding almost-dessins}
\label{findnear}
Suppose for example we want to find all (up to conjugation) planar 4-constellations  $[g_0,g_1,g_t,g_{\infty}]$ where the cycle structures of $g_1$ and $g_{\infty}$ are $(2^6)$ and $(3^4)$
($g_t$ is always a 2-cycle).
Up to conjugacy we may assume that $g_{\infty} = (1\,2\,3)(4\,5\,6)(7\,8\,9)(10\,\,11\,\,12)$.
The number of elements of $S_{12}$ of type $(2^6)$ is $10395$, and the number of 2-cycles is 66.
One could, for all $10395 \times 66$ combinations of $(g_1, g_t)$,
compute $g_0 = (g_1 g_t g_{\infty})^{-1}$, check if $[g_0,g_1,g_t,g_{\infty}]$ is transitive and planar, and if so, apply Algorithm 4.3. 
This works fine, but it can easily be sped up.

Since  $g_t$ is a 2-cycle and $<\hspace{-2pt}g_1,g_t,g_{\infty}\hspace{-2pt}>$ should be transitive,
it follows that $<\hspace{-2pt}g_1,g_{\infty}\hspace{-2pt}>$ may have at most two orbits in $\{1 \ldots 12\}$. So $g_1$ must connect some of the $g_{\infty}$-orbits $\{1,2,3\},\{4,5,6\},\{7,8,9\},\{10,11,12\}$.
Up to conjugation, we may assume $g_1$ connects the first two orbits with the 2-cycle $(1\,4)$
(we may also assume that $g_1$ contains either $(2\,7)$ or $(7\,\,10)$ since $g_1$ must connect more than one pair of $g_{\infty}$-orbits).
This way all almost-dessins (all 4-constellations up to conjugation) with such branching patterns can be found with little CPU time.

%

\subsection{Braid orbits}

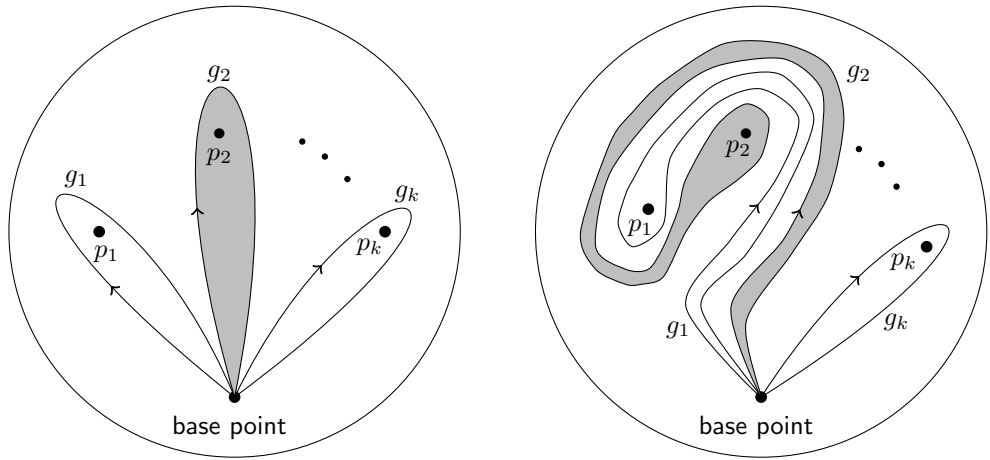
\begin{figure}[H]
\begin{tikzpicture}[smooth cycle]

\draw (1,2) circle (3cm);
\draw [fill] (1,-0.2) circle (2pt); \draw [fill] (-.8,2) circle (2pt);  
\draw [fill] (3,2) circle (2pt);
 \draw [fill] (1.9,3.2) circle (1pt); \draw [fill] (2.2,3) circle (1pt);\draw [fill] (2.5,2.7) circle (1pt);
\node at (0.93, -0.6){base point};
\draw[decoration={markings, mark=at position 0.3 with {\arrow[line width=.3mm]{>}}},postaction={decorate}](1,-.2) .. controls (-3.5,3.4) and (-.5,3.4) .. (1,-.2);
\draw[fill = lightgray] [decoration={markings, mark=at position 0.3 with {\arrow[line width=.3mm]{>}}},postaction={decorate}](1,-.2) .. controls (-.6,5.3) and (2.1,5.3) .. (1,-.2);
\draw (.8,3.3) node {$\bullet$};

 \draw [decoration={markings, mark=at position 0.3 with {\arrow[line width=.3mm]{>}}},postaction={decorate}](1,-.2) .. controls (2.7,3.3) and (5.3,3) .. (1,-.2);
 \node at (-1.1,2.7) {$g_1$}; \node at (.8,4.1) {$g_2$}; \node at (3.3, 2.5) {$g_k$};
 \node at (-.7, 1.7) {$p_1$}; \node at (.8,3) {$p_2$}; \node at (2.8, 1.8) {$p_k$};

\draw [fill] (8,-.2) circle (2pt);

\draw (8,2) circle (3cm);  \draw[fill] (6.5,2.3) circle (2pt);   \draw [fill] (10.2,1.8) circle (2pt);
\node at (7.93,-.6) {base point};

\draw [decoration={markings, mark=at position 0.1 with {\arrow[line width=.3mm]{>}}},postaction={decorate}]plot coordinates{(8,-.2) (7,1) (7.6,1.9) (8.2,2.7)(8.5,3.5) (8,3.9) (7.4,3.7) (6.8, 3) (6.7,2.1) (6.3,1.8) (6.1,2.1)(6.2,2.5) (6.3,2.7) (6.6,3.2)(7.3,3.9) (7.8,4.1) (8.3,4.07) (8.7,3.5)
                                                         (8.3,2.5) (7.8,1.7) (7.2,1)};
\draw[fill=lightgray][decoration={markings, mark=at position 0.06 with {\arrow[line width=.3mm]{>}}},postaction={decorate}] plot coordinates{(8,-.2) (7.6,1) (8.3,1.9) (8.7,2.7)(8.8,3.2) (8.8,3.9) (8.3,4.3) (7.4,4.2) (6.8, 3.8)  (6,2.8) (5.8,1.8) (6.2,1.5) (6.5,1.5) (6.7,1.8)(6.9,2.3)(7,2.7) (7.4, 3.4) (7.8,3.7) (8.1,3.5)(8,3) (7.7,2.7) (7.3,2.4)(7,2) (6.7, 1.4) (6.4, 1.3) (6.2, 1.3) (5.8, 1.5) (5.6,2) (6, 3.2) (6.3,3.6) (6.7,4) (7.5,4.5)(8.3,4.5) (8.9,4.2) (9.1,3.5)(8.9,2.5)(8.4,1.7)(7.8,.8)};
\draw (7.8,3.3) node {$\bullet$};

   \draw [decoration={markings, mark=at position 0.3 with {\arrow[line width=.3mm]{>}}},postaction={decorate}](8,-.2) .. controls (10.2,3.2) and (12.3,2.5) .. (8,-.2);
 \draw [fill] (9.3,3.1) circle (1pt); \draw [fill] (9.6,2.9) circle (1pt);\draw [fill] (9.8,2.6) circle (1pt);
 \node at (6.9,.7) {$g_1$}; \node at (9.3,4.1) {$g_2$}; \node at (9.8, .8) {$g_k$};
 \node at (6.4, 2.05) {$p_1$}; \node at (7.7,3.1) {$p_2$}; \node at (9.9, 1.6) {$p_k$};

\end{tikzpicture}
\caption{Action of $\sigma^2_1: [g_1 \ldots g_k] \mapsto [g_1^{\tau}, g_2^{\tau}, \, g_3 \ldots g_k]$  (here $\tau = g_1 g_2$).} \label{figure3}
\end{figure}

The braid group, generated by the braids $\sigma_1, \ldots, \sigma_{k-1}$, acts on $k$-constellations in the following way:
 \[ \sigma_i: [g_1 \ldots g_k] \mapsto [g_1 \ldots g_{i-1}, \,  g_{i+1},\, g^{-1}_{i+1}\,g_i\,g_{i+1}, \, g_{i+2} \ldots g_k]. \]
The points $p_i, p_{i+1}$ are swapped by $\sigma_i$ with a half-rotation. We will use
orbits under the {\em pure braid group} (Def. 9.11 in~\cite{RETBook}) which consists of products of $\sigma_i$'s that return $p_1,\ldots,p_n$  to their original locations.
The diagram in Figure~\ref{figure3}, taken from Section~1 in \cite{BraidActionPicture}, illustrates $\sigma^2_1$.
An algorithm is given in \cite{BraidOrbitComputation} for computing braid orbits of $k$-constellations $[g_1,\ldots,g_k]$.
Combining this with Algorithm~4.3 we obtain an algorithm that computes braid orbits of almost-dessins.

Our implementation \cite{dessindatabase} automatically generates Table~\ref{tableB1} below. 
First, it computes all planar \mbox{$\text{Belyi}^{(1)}$} branching patterns with $|\,E\,| = 5$ with the algorithm mentioned in Section 2.1.
Next, it computes all almost-dessins for these branching patterns as in Section~\ref{findnear}. {These are grouped into braid orbits labeled $N_1,\ldots,N_{68}$. 
Here 
$N_{66},N_{67},N_{68}$ are inside \mbox{$\text{Belyi}^{(2)}$} families.} 
The table shows the length (number of almost-dessins) of each orbit and its branching pattern.  Most branching patterns have one orbit, but some have zero or two.

\begin{table}[H]
{\small
  \begin{tabular}{c|l|l|c|l||c|l | c|c|l}
 $n$ & branching pattern & \hspace{-2pt}name\hspace{-2pt}  & $|\cal{O}|$ & \hspace{-2pt}decomp.\hspace{-3pt} & $n$ & branching pattern& \hspace{-2pt}name\hspace{-2pt}  & $|\cal{O}|$ & \hspace{-2pt}decomp. \\
 \hline
   &                                   &                 &      &            &    &                                   &          &      &\\[-8pt]
 2  & $(2)$, $(1^2)$, $(1^2)$           & $N_{1}$         & 1    &            & 8  & $(1^2,2,4)$, $(2^4)$, $(2,3^2)$   & $N_{38}$ & 9    &\\
   & $(1^2)$, $(1^2)$, $(2)$           & $N_{2}$         & 1    &            &    & $(1^2,3^2)$, $(2^4)$, $(2,3^2)$   & $N_{39}$ & 4    & $4\circ 2$\\
3  & $(1,2)$, $(1^3)$, $(3)$           & $N_{3}$         & 1    &            &    & $(1,2^2,3)$, $(2^4)$, $(2,3^2)$   & $N_{40}$ & 9    &\\
   & $(3)$, $(1,2)$, $(1^3)$           & $N_{4}$         & 1    &            &    & $(2^4)$, $(2^4)$, $(2,3^2)$       & ---      &  --- &\\
   & $(1,2)$, $(1,2)$, $(1,2)$         & $N_{5}$         & 4    &            & 9  & $(1^3,6)$, $(1,2^4)$, $(3^3)$     & $N_{41}$ & 3    & $3\circ 3$\\
4  & $(1,3)$, $(1^2,2)$, $(1,3)$       & $N_{6}$         & 6    &            &    &  \hspace{2pt}idem                 & $N_{42}$ & 9    &\\
   & $(2^2)$, $(1^2,2)$, $(1,3)$       & $N_{7}$         & 3    &            &    & $(1^2,2,5)$, $(1,2^4)$, $(3^3)$   & $N_{43}$ & 18   &\\
   & $(4)$, $(2^2)$, $(1^4)$           & $N_{8}$         & 1    & $2\circ 2$ &    & $(1^2,3,4)$, $(1,2^4)$, $(3^3)$   & $N_{44}$ & 15   &\\
   & $(1,3)$, $(2^2)$, $(1^2,2)$       & $N_{9}$         & 3    &            &    & $(1,2^2,4)$, $(1,2^4)$, $(3^3)$   & $N_{45}$ & 12   & $3\circ 3$\\
   & $(2^2)$, $(2^2)$, $(1^2,2)$       & $N_{10}$        & 2    & $2\circ 2$ &    & $(1,2,3^2)$, $(1,2^4)$, $(3^3)$   & $N_{46}$ & 12   &\\
   & $(1^2,2)$, $(2^2)$, $(2^2)$       & $N_{11}$        & 2    & $2\circ 2$ &    & $(2^3,3)$, $(1,2^4)$, $(3^3)$     & $N_{47}$ & 3    & $3\circ 3$\\
   & $(1^4)$, $(2^2)$, $(4)$           & $N_{12}$        & 1    & $2\circ 2$ & 10 & $(1^3,7)$, $(2^5)$, $(1,3^3)$     & $N_{48}$ & 15   &\\
5  & $(1,4)$, $(1,2^2)$, $(1^2,3)$     & $N_{13}$        & 10   &            &    & $(1^2,2,6)$, $(2^5)$, $(1,3^3)$   & $N_{49}$ & 15   &\\
   & $(2,3)$, $(1,2^2)$, $(1^2,3)$     & $N_{14}$        & 7    &            &    & $(1^2,3,5)$, $(2^5)$, $(1,3^3)$   & $N_{50}$ & 15   &\\
   & $(1^2,3)$, $(1,2^2)$, $(2,3)$     & $N_{15}$        & 7    &            &    & $(1^2,4^2)$, $(2^5)$, $(1,3^3)$   & $N_{51}$ & 12   &\\
   & $(1,2^2)$, $(1,2^2)$, $(2,3)$     & $N_{16}$        & 6    &            &    & $(1,2^2,5)$, $(2^5)$, $(1,3^3)$   & $N_{52}$ & 15   &\\
6  & $(1^2,4)$, $(1^2,2^2)$, $(3^2)$   & $N_{17}$        & 3    & $3\circ 2$ &    & $(1,2,3,4)$,\,$(2^5)$,\,$(1,3^3)$ & $N_{53}$ & 18   &\\
   &  \hspace{2pt}idem                 & $N_{18}$        & 6    &            &    & $(1,3^3)$, $(2^5)$, $(1,3^3)$     & ---      & ---  &\\
   & $(1,2,3)$, $(1^2,2^2)$, $(3^2)$   & $N_{19}$        & 12   &            &    & $(2^3,4)$, $(2^5)$, $(1,3^3)$     & ---      & ---  &\\
   & $(2^3)$, $(1^2,2^2)$, $(3^2)$     & $N_{20}$        & 3    & $3\circ 2$ &    & $(2^2,3^2)$, $(2^5)$, $(1,3^3)$   & $N_{54}$ & 6    &\\
   & $(1,5)$, $(2^3)$, $(1^3,3)$       & $N_{21}$        & 5    &            & 12 & $(1^4,8)$, $(2^6)$, $(3^4)$       & $N_{55}$ & 3    & $3 \circ 2 \circ 2$\\
   & $(2,4)$, $(2^3)$, $(1^3,3)$       & $N_{22}$        & 2    & $2\circ 3$ &    &     \hspace{2pt}idem              & $N_{56}$ & 4    &\\
   & $(3^2)$, $(2^3)$, $(1^3,3)$       & $N_{23}$        & 2    &            &    & $(1^3,2,7)$, $(2^6)$, $(3^4)$     & $N_{57}$ & 7    &\\
   & $(1^2,4)$, $(2^3)$, $(1,2,3)$     & $N_{24}$        & 9    &            &    & $(1^3,3,6)$, $(2^6)$, $(3^4)$     & $N_{58}$ & 4    & $4\circ 3$\\
   & $(1,2,3)$, $(2^3)$, $(1,2,3)$     & $N_{25}$        & 10   &            &    & $(1^3,4,5)$, $(2^6)$, $(3^4)$     & $N_{59}$ & 10   &\\
   & $(2^3)$, $(2^3)$, $(1,2,3)$       & $N_{26}$        & 2    & $2\circ 3$ &    & $(1^2,2^2,6)$, $(2^6)$, $(3^4)$   & $N_{60}$ & 9    & $3\circ 4$\\
7  & $(1^2,5)$, $(1,2^3)$, $(1,3^2)$   & $N_{27}$        & 21   &            &    & $(1^2,2,3,5)$, $(2^6)$, $(3^4)$   & $N_{61}$ & 15   &\\
   & $(1,2,4)$,\,$(1,2^3)$,\,$(1,3^2)$ & $N_{28}$        & 24   &            &    & $(1^2,2,4^2)$, $(2^6)$, $(3^4)$   & $N_{62}$ & 6    & $3 \circ 2 \circ 2$\\
   & $(1,3^2)$, $(1,2^3)$, $(1,3^2)$   & $N_{29}$        & 12   &            &    & $(1^2,3^2,4)$, $(2^6)$, $(3^4)$   & ---      & ---  &\\
   & $(2^2,3)$, $(1,2^3)$, $(1,3^2)$   & $N_{30}$        & 9    &            &    & $(1,2^3,5)$, $(2^6)$, $(3^4)$     & ---      & ---  &\\
8  & $(1^2,6)$, $(2^4)$, $(1^2,3^2)$   & $N_{31}$        & 4    & $4\circ 2$ &    & $(1,2^2,3,4)$, $(2^6)$, $(3^4)$   & $N_{63}$ & 9    & $3\circ 4$\\
   &   \hspace{2pt}idem                & $N_{32}$        & 12   &            &    & $(1,2,3^3)$, $(2^6)$, $(3^4)$     & $N_{64}$ & 4    & $4\circ 3$\\
   & $(1,2,5)$, $(2^4)$, $(1^2,3^2)$   & $N_{33}$        & 10   &            &    & $(2^4,4)$, $(2^6)$, $(3^4)$       & $N_{65}$ & 3    & $S_3 \circ 2$\\
   & $(1,3,4)$, $(2^4)$, $(1^2,3^2)$   & $N_{34}$        & 15   &            &    & $(2^3,3^2)$, $(2^6)$, $(3^4)$     & ---      & ---  &\\
   & $(2^2,4)$, $(2^4)$, $(1^2,3^2)$   & $N_{35}$        & 6    & $2\circ 4$ & 4  & $(1^4)$, $(4)$, $(1,3)$     & $N_{66}$ & 1    &\\
   & $(2,3^2)$, $(2^4)$, $(1^2,3^2)$   & $N_{36}$        & 4    & $4\circ 2$ & 6  & $(1^4,2)$, $(2,4)$, $(3^2)$     & $N_{67}$ & 4    &\\
   & $(1^3,5)$, $(2^4)$, $(2,3^2)$     & $N_{37}$        & 6    &            &    & $(1^4,2)$, $(2^3)$, $(6)$       & $N_{68}$ & 2    & $2 \circ 3$ \\

\end{tabular}
}
\caption{Braid orbits of the almost-dessins with $|\,E\,| = 5$} \label{tableB1}
\end{table}

Column $|\cal{O}|$ in the above table gives the number of almost-dessins in each braid orbit. The complete table of almost-dessins themselves is given in \cite{dessindatabase}. The notation $3 \circ 4$ means that any $\text{Belyi}^{(1)}$ map $f$ for this orbit equals $g \circ h$ for some $g,h$ of degrees 3,4.
The notation $S_3 \circ 2$ means $f = g \circ h$ where $g$ has three $3 \circ 2$-decompositions and one $2 \circ 3$-decomposition.
We do not need explicit $f \in \mathbb{C}(x)$ in order to find any of the information listed in Table~\ref{tableB1},
including the decomposition structure of $f$ (the almost-dessins $[g_0,g_1,g_t,g_{\infty}]$ suffice).
Decompositions of $f$ correspond to subfields $\mathbb{C}(f) \subseteq E \subseteq \mathbb{C}(x)$, which in turn
correspond to subgroups of $G := <\hspace{-2pt}g_0,g_1,g_t,g_{\infty}\hspace{-2pt}>$ that contain ${\rm Stab}(1) = \{g \in G | g(1) = 1\}$.

\pagebreak

The sections below can use Table~\ref{tableB1} to prove that our database covers all $\text{Belyi}^{(1)}$ maps,
as everything in Table~\ref{tableB1} was computed independently of these functions.

\subsection{Continuation of Example 2.1} \label{section5.3}
Let $B$ be the third branching pattern under $n=10$ in Table~\ref{tableB1}.
Example~2.1 gave a $\text{Belyi}^{(1)}$ map $f(s,x)$
for $B$.
The table shows that $B$ has 15 distinct almost-dessins, in one braid orbit named $N_{50}$.
Let $\phi_f(s) \in \mathbb{Q}(s)$ be the rational function of degree 15 in Equation~(\ref{valuet}), as in Definition~\ref{specialvalues}.
Choose any $t_0 \in \mathbb{P} - \{0,1,\infty\}$, and let $S := \phi_f^{-1}(\{t_0\})   \subset \mathbb{P}^1$.
If $\alpha \in S$, then $f(\alpha,x)$ is a $\text{Belyi}^{(1)}$ map for $B$ 
that ramifies only above $\{0,1,t_0,\infty\}$,
{\em assuming $f$'s family is gap-free} as in Definition~\ref{specialvalues}.
One could compute ($2 \mapsto 1$ in Section~\ref{dessins})
the almost-dessin $D_{\alpha}$ of $f({\alpha},x)$ for each $\alpha \in S$, then take $D_* := \{D_{\alpha} | \alpha \in S\}$, and check that $N_{50} = D_*$.
However, it is not hard to see that this check is not necessary for this $B$.


Let $\gamma$ be a loop in $\mathbb{P} - \{0,1,\infty\}$ with base point $t_0$.
Applying analytic continuation to $\phi_f^{-1}(\{t\})$, with $t$ following $\gamma$, 
gives a map from $S$ to $S$.
This gives an action of the fundamental group $\pi_1(\mathbb{P} - \{0,1,\infty\}, t_0)$ on $S$.
Since $D_*$ is an image of $S$, the fundamental group acts on $D_*$ as well.
Figure~\ref{figure3} illustrates how this corresponds to an action of the pure braid group.

Table~\ref{tableB1} implies $D_* \subseteq N_{50}$ because according to Table~\ref{tableB1},
all 15 almost-dessins for $B$ are in $N_{50}$.  
%
Then $D_* = N_{50}$ because the pure braid group acts on $D_*$ and $N_{50}$ is an orbit.

\begin{proposition} \label{PropFamily} Let $f$ and $B$ be as above.  If $g \in \mathbb{C}(x)$ has branching pattern $B$ then it is M\"obius-equivalent to a member of $f$'s family.
\end{proposition}
\noindent {\bf Proof:} Let $t_0$ be the branch point of $g$ not in $\{0,1,\infty\}$. $D_* = N_{50}$ (we checked that $f$'s family is gap-free \cite{dessindatabase}).
The almost-dessin of $g$ has branching pattern $B$, is thus in $N_{50}$ and hence equals $D_{\alpha}$ for some $\alpha \in S$. Then $f(\alpha,x)$ is M\"obius-equivalent to $g$,
see correspondence~(\ref{fcor}) in Section~3. \eindebewijs




\subsection{A branching pattern with two orbits} \label{twobelyi-1mapswithsamepattern}
Let
$$f_1 = \frac{3\,s(x-1)x^2+4}{4\,(s(x-1)x^2+1)^3} \ \ {\rm and} \ \
    f_2= -\frac{s^2\,((4s-3)x^3+6\,(s-1)x^2+3\,(3s^2-2s-1)x-4s))}{4\,(x^3+2x^2+(2s+1)\,x+s)^3}.    $$
Both are gap-free, have branching $(1^3, 6), (1,2^4), (3^3)$ above $0, 1, \infty$
and one more branch point $t = \phi_{f_1}(s)$ and $t = \phi_{f_2}(s)$ respectively.
The degree of $\phi_{f_1}$ is 3. This, combined with argument from Section~\ref{section5.3}, suffices to prove that $f_1$ covers $N_{41}$ in Table~\ref{tableB1}.
However, the fact that $\phi_{f_2}$ has degree 9 is not enough to demonstrate that $f_2$ covers $N_{42}$ 
because, in the notations from Section~\ref{section5.3}, the cardinality of $\{ D_{\alpha} | \alpha \in S \}$ could be less than the cardinality of $S$.

%
%
%

\begin{defn} \label{dupfree}
A $\text{Belyi}^{(1)}$ map $f \in \mathbb{C}(s)(x)$ is called {\em duplicate-free} if $\{ f(\alpha, x)\,|\, \alpha \in \phi_f^{-1}( \{t_0\} ) \}$
has ${\rm deg}_s(\phi_f)$ distinct almost-dessins for any $t_0 \not\in \{0,1,\infty\}$.
\end{defn}
After verifying that $f_2$ is duplicate-free we may conclude that it covers $N_{42}$, since it is the only orbit for this branching pattern of length $9$.
\begin{remark}
If $f \in \mathbb{C}(s)(x)$ is a duplicate-free $\text{Belyi}^{(1)}$ map then $\phi_f \in \mathbb{C}(s)$ is a Belyi map, and its dessin can be computed directly from a 4-constellation $[g_0,g_1,g_t,g_{\infty}]$ of $f$.
\end{remark}
{\bf Proof:} Definition~\ref{dupfree} immediately implies $ | \phi_f^{-1}( \{t_0\} ) | \geq {\rm deg}_s( \phi_f )$ for any $t_0 \not\in \{0,1,\infty\}$, in other words, $\phi_f$ is a Belyi map.
Take braid actions that correspond to looping $t$ around $0,1,\infty$. Let $h_0,h_1,h_{\infty}$ be the corresponding permutations of the almost-dessins, then $[h_0,h_1,h_{\infty}]_\sim$ is the dessin of $\phi_f$.
\eindebewijs
It was fortunate these dessins were always planar in our database, otherwise our $\text{Belyi}^{(1)}$ maps could not have been in $\mathbb{Q}(s)(x)$, complicating the algorithms.


\subsection{Proving completeness of our table of $\text{Belyi}^{(1)}$ maps}
To our surprise, Section~\ref{section2} often produced $\text{Belyi}^{(1)}$ maps $f$ that were not duplicate-free, where
the degree of $\phi_f$ was twice the number of distinct almost-dessins. For such cases, we computed
automorphisms $\tau \in {\rm Aut}(\mathbb{Q}(s))$ of order 2 for which $\tau(\phi_f) = \phi_f$, in order to find $\tau$
for which $\tau(f)$ is M\"obius-equivalent to $f$. Let $\tilde{s}$ be a generator of the subfield of $\mathbb{Q}(s)$ fixed by $\tau$.

We write $\phi_f$ as element of $\mathbb{Q}(\tilde{s})$ and use it to
search for a $\tilde{f}(s,x)$ for which $\tilde{f}(\tilde{s},x)$ is M\"obius-equivalent to $f$.
Then $\phi_{\tilde{f}}$ has half the degree of $\phi_f$. This way, we managed to make every member of our $\text{Belyi}^{(1)}$ table duplicate-free.
%
After suitable M\"obius transformations, we managed to make them gap-free as well.
The arguments of the previous two subsections now suffice to prove that our $\text{Belyi}^{(1)}$ table \cite{dessindatabase} is complete. But we implemented a more direct verification as well:

Let $F^{(1)}_1 \ldots F^{(1)}_{68}$ be the explicit $\text{Belyi}^{(1)}$ maps at \cite{dessindatabase}. For each $i$ we check
that $F^{(1)}_i$ is gap-free, compute almost-dessin for $F^{(1)}_i$ and check that it is in $N_i$. This suffices to prove that, up to M\"obius-equivalence,
the families of $F^{(1)}_1 \ldots F^{(1)}_{68}$ contain all rational $\text{Belyi}^{(1)}$ maps with $|\,E\,| = 5$.
We also compute the degree of $\phi_{F^{(1)}_i}$ and check that it equals {$|\cal{O}|$ which denotes the number of elements of the braid orbit.

{As for the Belyi case, we also implemented a faster approach, based on five point invariants, to prove the completeness of the $\text{Belyi}^{(1)}$ tables.
Here we used not one, but two algebraically independent five point invariants $I_5$ and $\tilde{I}_5$. Section \ref{FindF} gives more details about these invariants.
For each $f$ in the $\text{Belyi}^{(1)}$ table, $I_5(f)$ and $\tilde{I}_5(f)$ are both in $\mathbb{Q}(s)$, and thus satisfy an algebraic relation. 
If two $\text{Belyi}^{(1)}$ maps give distinct algebraic relations, then they can not be part of the same family.
This turned out to be the case for any pair in our table with the same branching pattern.}

\subsection{Belyi maps inside $\text{Belyi}^{(1)}$ families} \label{BB1}
The family of $F^{(1)}_i$ contains, up to M\"obius-equivalence, all $\text{Belyi}^{(1)}$ maps with almost-dessin in $N_i$.
But it often contains Belyi maps as well; if $s_0$ is special (Definition~\ref{specialvalues}) then $F^{(1)}_i(s_0,x)$ is a Belyi map.
Depending on whether $\phi_{F^{(1)}_i}(s_0)$ is 0, 1, or $\infty$, the dessin of this Belyi map is $[g_0 g_t^\tau,\,  g_1, g_{\infty}]_\sim$ (where $\tau = g_1^{-1}$),
$[g_0,\, g_1g_t,\, g_{\infty}]_\sim$ or $[g_0, g_1,\, g_tg_{\infty}]_\sim$.
{All dessins for which $|E|$ is larger in the definition from Remark 1.4 than in our definition 1.2
can be obtained this way. So one would expect that the families of $F^{(1)}_1 \ldots F^{(1)}_{68}$ contain
a Belyi maps for each of those dessins.
However, there are a few exceptions; some dessins that can be obtained this way from $N_i$ do not {\em directly} appear in $F^{(1)}_i$'s family because they correspond to degenerate values of $s$.
They do appear {\em indirectly}, i.e. in another, less favorable, parametrization of $F^{(1)}_i$:
\begin{exmp} \label{IndirectBelyi}
Let $f := -27(s x^4-2 s x^3+s x^2+1)^2 /(s x^4-2 s x^3+s x^2-3)^3$. 
Let $\tilde{s} := \sqrt[4]{s}$ and \linebreak $\tilde{f} := f(s, x/\tilde{s})$.
The point $\tilde{s}=0$ is degenerate for $f$ but special for $\tilde{f}$, where it evaluates to a Belyi map $g = -27(x^4+1)^2 /(x^4-3)^3$.
Although $f$'s family is a proper subset of $\tilde{f}$'s family,
we prefer $f$ 
because it is duplicate-free.
\end{exmp}

\section{$\text{Belyi}^{(2)}$ maps} \label{belyi-2maps}
\begin{defn}
Let $S$ be a subset of $\mathbb{P}^1$ with $n$ elements. The $n$-{\em point-polynomial} 
$P_S \in \mathbb{C}[x]$ is the product of $x-p$ taken over all $p \in S-\{\infty\}$.
It has degree $n-1$ if $\infty \in S$ and degree $n$ otherwise.
Let $k$ be a subfield of $\mathbb{C}$. We say $S$ is {\em defined over} $k$ if $P_S \in k[x]$.
\end{defn}
\noindent If $f \in k(x)$ then its set of $(k,\ell,m)$-exceptional points is defined over $k$.
Let
\[ F_4^{(2)}(a,b,c,d,x) = 1-\frac{(x^2+ax+b)^2}{c(x+d)^3}, \ \  F_6^{(2)}(a,b,c,d,x) :=
1-\frac{(x^3 + 3 a x^2+ bx +c)^2}{(x^2+2ax+d)^3}. \]
Their branching patterns are $B_4 = (1^4), (2^2), (1,3)$ and $B_6 = (1^4,2), (2^3), (3^2)$.
Both are two-dimensional families up to M\"obius-equivalence (two of the 4 parameters $a,b,c,d$ can be eliminated
with a linear transformation on $x$).

\begin{lemma} \label{B2l} Let $k$ be a subfield of $\mathbb{C}$ and $f \in k(x)$ a $\text{Belyi}^{(2)}$ map with $|\,E\,| = 5$.
Then $f$ has branching pattern $B_4$ or $B_6$, and there exist unique
$m \in \{1/(x-p) \,|\, p \in k\} \bigcup \{x\}$ and $a,b,c,d \in k$ such that $f$ equals $F_4^{(2)}(a,b,c,d, m)$ if $f$ has $B_4$,
and $F_6^{(2)}(a,b,c,d, m)$ if $f$ has $B_6$.
\end{lemma}
\noindent {\bf Proof:} Our implementation mentioned in Section~2.1 shows (it is also easy to show by hand) that $B_4$ and $B_6$ are the only planar $\text{Belyi}^{(2)}$
branching patterns with $|\,E\,| = 5$.
If $f$ has $B_6$, then $(1^4,2)$ indicates that it has a unique root $p$ of order 2, and four roots of order 1.
The part $(2^3)$ of $B_6$ indicates that numerator of $1-f$ must be a square, while  $(3^2)$ indicates that the denominator is a cube.
If $p = \infty$ then $f(\infty)=0$ with multiplicity 2, which implies that the numerator and denominator of $1-f$ must have the same degree, same leading
coefficient, and the same $x^5$-coefficient as well. Then $f$ must
equal $F_6^{(2)}(a,b,c,d,x)$ for some $a,b,c,d$, uniquely determined by $f$, and hence in $k$.
If $p \neq \infty$, the M\"obius-transformation $m$ moves $p$ to $\infty$, after which the same argument applies.

If $f$ has $B_4$, then let $p$ be the unique pole of order 1. If $p=\infty$, then
the denominator of $f$ must be a cube and the numerator of $1-f$ a square, hence $f = F_4^{(2)}(a,b,c,d,x)$ for unique $a,b,c,d \in k$.
The case $p \neq \infty$ again reduces to this under $m$. \eindebewijs

As there are only two cases, it is not hard to
solve the $\text{Belyi}^{(2)}$ part of goal~(c) from the introduction:

\begin{framed}
\noindent \textbf{Algorithm 6.1: FindBelyi2} \\
\textbf{Input:} A field $k \subseteq \mathbb{C}$ and a 5-element subset $S = \{q_1 \ldots q_5\} \subset \mathbb{P}^1$ defined over $k$. \\
\textbf{Output includes:} Every $\text{Belyi}^{(2)}$ $f \in k(x)$ such that $E(f) = S$. \\[3pt]
{For each} $p$ in $S  \bigcap (k \bigcup \{\infty\})$ do: \\[3pt]
\mbox{} \ \ \ \textbf{Step 1}. Let $m$ be as in Lemma~\ref{B2l} and $\tilde{m}$ be its inverse ($x$ if $p=\infty$, otherwise $1/x + p$). \\
\mbox{} \ \ \ \textbf{Step 2}. Comparing the numerator of $F_4^{(2)}(a,b,c,d,x)$ with $m(S)$ gives 4 equations in $a,b,c,d$. \\
\mbox{} \ \ \ \textbf{Step 3}. Two equations are linear in a variable, solving these leaves 2 equations in 2 unknowns. \\
\mbox{} \ \ \ \textbf{Step 4}. Compute all solutions over $k$ with a resultant. \\
\mbox{} \ \ \ \textbf{Step 5}. For each solution, append $F_4^{(2)}(a,b,c,d, m)$ to the output. \\
\mbox{} \ \ \ \textbf{Step 6}. Doing the same for $F_6^{(2)}(a,b,c,d,x)$ gives 4 equations, one of which is linear. \\
\mbox{} \ \ \ \textbf{Step 7}. With a pre-computed \cite{VJKThesis} elimination we obtain an equation of degree 12 for $a$. \\
\mbox{} \ \ \ \textbf{Step 8}. After computing its roots in $k$, two equations in two unknowns remain. \\
\mbox{} \ \ \ \textbf{Step 9}. Compute solutions as in Step~4 and for each, append $F_6^{(2)}(a,b,c,d, m)$ to the output.
\end{framed}
The program finds all $\text{Belyi}^{(2)}$ maps for $S$ in $k(x)$
but it also finds certain Belyi or $\text{Belyi}^{(1)}$ maps:
$F^{(1)}_{66}$ is a special case of $F^{(2)}_4$ while $F^{(1)}_{67}$ and $F^{(1)}_{68}$ are special cases of $F^{(2)}_6$.
							%
							%
So we can remove these three from our $\text{Belyi}^{(1)}$ table without interfering with goal~(c). 
To cover goal~(c) for Belyi and $\text{Belyi}^{(1)}$ maps we need one more ingredient, which will be the topic of the next section.

\section{Five point invariants}
\label{FindF}
Given $k$ and $S$, our goal is to quickly find, if it exists, $f \in k(x)$ such that $E(f) = S$.
After running Algorithm FindBelyi2 we may assume that $f$ is Belyi or $\text{Belyi}^{(1)}$.
Such $f$ must be M\"obius-equivalent to a member of our Belyi or $\text{Belyi}^{(1)}$ table because they were proved to be complete.

It is not efficient to search for a M\"obius-equivalence between $S$ and the exceptional points
of each of the many entries of the Belyi table. For the $\text{Belyi}^{(1)}$ table, one first needs to
find the correct value of the parameter $s$ before a M\"obius-equivalence could occur.

Let $k_5$ be the set of 5-element subsets $S \subset \mathbb{P}^1$ that are defined over $k$.
A {\em five-point-invariant} is a function $k_5 \rightarrow k$ that is invariant under M\"obius-transformations. We implemented two such functions.
The first, called $I_5$, maps $S$  to $\sum_{q \in S} j(S - \{q\})$ where $j(T)$ refers to the j-invariant
of a set $T$ with $4$ points. More precisely, if $T = \{q_1,q_2,q_3,q_4\}$ then $j(T)$ is the j-invariant of the elliptic curve $y^2 = (x-q_1)(x-q_2)(x-q_3)(x-q_4)$,
where a factor $x-q_i$ is omitted if $q_i = \infty$.
The second invariant $\tilde{I}_5$ is similar, except that it uses the sum of the squares of the $j$-invariants.

If $f$ has $5$ exceptional points $S = \{q_1,\ldots,q_5\}$, then $I_5(f)$ denotes $I_5( S )$.
We attach $I_5(f)$ to  each Belyi map $f$ in our database.
To each $\text{Belyi}^{(1)}$ map $f \in \mathbb{Q}(s)(x)$, we attach $I_5(f)$ and $\tilde{I}_5(f)$, which are elements of $\mathbb{Q}(s)$.
For a Belyi map $f$, the invariant $I_5(f)$ is either a rational or an algebraic number (we insert its
minimal polynomial over $\mathbb{Q}$ into the table). We do not use five-point invariants for $\text{Belyi}^{(2)}$ maps because there were only two cases.

These invariants give an efficient solution to goal (c), they rapidly eliminate nearly all entries that do not lead to a solution.

\begin{framed}
\noindent \textbf{Algorithm 7.1: FindF (goal~(c))} \\
\textbf{Input:} A field $k \subseteq \mathbb{C}$ and a $5$-element subset $S = \{q_1 \ldots q_5\} \subset \mathbb{P}^1$ defined over $k$. \\
\textbf{Output:} Every element of $f \in k(x)$ such that $E(f) = S$. \\[3pt]
\textbf{Step 1}. $A := $\,FindBelyi2($S) \subset k(x)$, \ \  $i_5 :=  I_5(S) \in k$, \ \  $\tilde{i}_5 :=  \tilde{I}_5(S) \in k$.  \\
\textbf{Step 2}. For each $f$ in the Belyi table whose $I_5$ matches $i_5$, adjoin $f(m)$ to $A$ for every (if any) \\
\mbox{} \hspace{32pt} M\"obius-transformation $m$ that sends $S$ to $E(f)$. \\
\textbf{Step 3}. For each $f$ in the $\text{Belyi}^{(1)}$ table, compute the gcd of the numerators of $I_5(f) - i_5$  \\
\mbox{} \hspace{32pt} and $\tilde{I}_5(f) - \tilde{i}_5$. If this gcd is not 1, then compute all its roots in $k$. For each \\
\mbox{} \hspace{32pt} non-degenerate root $s_0$, evaluate $f$ at $s = s_0$ and then proceed as in Step 2. \\
\textbf{Step 4}. Return $A$.
\end{framed}
For each $f$ in our Belyi table, if $\alpha = I_5(f)$, then $f$ turned out to be in $\mathbb{Q}(\alpha)(x)$. But if for example $\alpha = {\rm RootOf}(x^2-x-1)$ while the input of {\tt FindF} is defined over say $k = \mathbb{Q}(\sqrt{5})$,
then we must replace $\alpha$ by its corresponding element(s) of $k$ before one can use $f$ (use $\alpha \mapsto i_5$ to map $f$ to an element of $k(x)$).
Computing $m$ requires some care too,  for details see our implementation \cite{dessindatabase}.
In Step~3 it is important that every member of our $\text{Belyi}^{(1)}$ table is duplicate-free, this ensures that if a $\text{Belyi}^{(1)}$ map in $k(x)$ has $5$ exceptional points, then the corresponding value of $s$ is unique and thus in $k \bigcup \{\infty\}$.
The algorithm does not consider $s=\infty$ since it is degenerate for every member of our $\text{Belyi}^{(1)}$ table.



\section{goal (d)} \label{motivation}



The Gauss Hypergeometric Function $\mbox{}_2F_1(a,b ; c\,|\,x)$ satisfies the so-called Gauss Hypergeometric Equation
\begin{equation}
	x(1-x)y'' + (c-(a+b+1)x)y'-aby = 0.
	\label{GHE}
\end{equation}
It has singularities at $0, 1, \infty$ with {\em exponents} $\{0,1-c\}$, $\{0,c-a-b\}$, $\{a,b\}$ respectively. The exponent differences are $(e_0,e_1,e_\infty) = (1-c, c-a-b, b-a)$ up to sign. The numbers $(k,\ell,m)$ from Definition~\ref{def11} correspond to a GHE (equation~(\ref{GHE})) with the following exponent differences:
\begin{equation}
	(e_0,e_1,e_\infty) = (1/k,1/\ell,1/m).
	\label{RelationExpDiff}
\end{equation}

\pagebreak 

Finding a $_2F_1$-type solution of a second order differential equation $L$ is equivalent to  finding a combination of transformations~(i),(ii),(iii) that sends
the GHE~(\ref{GHE})  to $L$:\\

  \hspace{2pt} (i) Change of variables: $y(x)\, \mapsto\, y(f)$

  \hspace{2pt} (ii) Gauge transformation:  $y \,\mapsto \, r_0 y + r_1 y'$

  \hspace{2pt} (iii) Exponential product: $y \,\mapsto \,\exp(\int r \,dx)\cdot y$  \ (in Conjecture~1, $\exp(\int r \,dx)$ will be algebraic). \\


  Let $L$ be as in Conjecture~1, with coefficients $a_i \in k(x)$ for some field $k \subseteq \mathbb{C}$, and with $5$ true singularities $S = \{q_1 \ldots q_5\}$, at least one of them logarithmic. Our tasks are (1): Use algorithm {\tt FindF} to find
  $f \in k(x)$ (if it exists) such that $E(f) = S$ and (2): Find a combination of transformations (i),(ii),(iii) that sends the GHE~(\ref{GHE})  with $(e_0,e_1,e_\infty) = (0,\frac12,\frac13)$  to $L$.

\subsection{Example}

Let $L$ be:
$$ y'' + {\frac { \left( 8\,{x}^{4}-{x}^{2}+2\,x-3 \right) }{x \left( x+1 \right)  \left( 4\,x+3 \right)  \left( {x}^{2}-2\,x+
3 \right) }}y' - {\frac {4\,{x}^{2}}{ \left( {x}^{2}-2\,x+3 \right) ^{2}
 \left( x+1 \right) ^{2} \left( 4\,x+3 \right) }}y = 0
$$



\begin{enumerate}
\item 
Find the true (= non-removable) singularities \cite{VJKThesis}. In this example, all singularities except $x = \infty$.  \\
The $5$-point polynomial is $P = x(x+1)(x+3/4)(x^2-2x+3)$.



\item 
{\tt FindF} finds the following functions $f$ such that $E(f)$ are given by $P$:

$F_{\rm list}:= [ {\frac {-{x}^{8}}{4 \left( {x}^{2}-2\,x+3 \right)
 \left( 4\,x+3 \right)  \left( x+1 \right) ^{2}}},{\frac { -4\left(
x+1 \right) ^{2} \left( {x}^{2}-2\,x+3 \right) {x}^{4}}{ \left( 4\,x+3
 \right) ^{2}}},{\frac { \left( x+1 \right) ^{4} \left( {x}^{2}-2
\,x+3 \right) ^{2}}{ 4\left( 4\,x+3 \right) {x}^{4}}},{\frac {
 -64\left( {x}^{2}-2\,x+3 \right)  \left( x+1 \right) ^{2}{x}^{12}}{
 \left( 4\,x+3 \right)  \left( 8\,{x}^{4}+36\,x+27 \right) ^{3}}},$

${\frac {64 \left( x+1 \right) ^{6} \left( {x}^{2}-2\,x+3 \right) ^{3}{x}
^{4}}{ \left( 4\,x+3 \right)  \left( 8\,{x}^{4}-4\,x-3 \right) ^{3}}},
1+3\,{\frac { \left( {x}^{2}-10\,x-3 \right) ^{2}}{ \left( 5\,x-3
 \right) ^{3} \left( x+1 \right) }}]$.


\item 
$H: \, {y''}-{\frac { \left( -3/2\,x+1 \right)}{x \left( x-
1 \right) }}{y'}+{\frac {5}{144}}\,{\frac {1}{x \left( x-1 \right) }}y = 0$ \, is the GHE~(\ref{GHE})  with $(e_0,e_1,e_\infty) = (0,\frac12,\frac13)$.
Among its solutions is $y(x)=\mbox{}_2F_1\Big(\frac{1}{12},\frac{5}{12}; 1 \,|\, x \Big)$.
We need to find transformations that send $H$ to $L$. For each $f \in F_{\rm list}$, apply change of variables $x \mapsto f$ to $H$, and then find
transformations (ii)+(iii) \cite{VJKThesis}.

\item 
Transformations (ii)+(iii) only exist for the third element in $F_{\rm list}$. Applying transformation (i), \newline $x \mapsto f = {\frac { \left( x+1 \right) ^{4} \left( {x}^{2}-2\,x+3 \right) ^{
2}}{4 \left( 4\,x+3 \right) {x}^{4}}}$, \ to $H$ produces \ $H_f: \ y'' + {\frac { ( 10\,{x}^{4}-{x}^{2}-6\,x-9 )}{x ( 4\,x+3 )  ( {x}^{2}-2\,x+3 )
 ( x+1 ) }} y' +5\,{\frac {( x+1) ^{2}}{{x}^{2}( 4\,x+3) ^{2}}}y = 0.$


\item 
$y(f)$ is a solution
of $H_f$.
Computing transformations (ii)+(iii) gives a solution of $L$:\\
$Y = { \frac{(\frac{x+1}{x})^{1/3}(x^2-2x+3)^{1/6}}{(4x+3)^{1/12}}} \cdot ~_2F_1\Big(\frac{1}{12},\frac{5}{12}; 1 \,|\,{\frac { \left( x+1 \right) ^{4} \left( {x}^{2}-2\,x+3 \right) ^{
2}}{4 \left( 4\,x+3 \right) {x}^{4}}} \Big)$.  To obtain another solution, replace $y(x)$ by another solution of $H$.
\end{enumerate}
Our implementation \cite{dessindatabase} performs the above steps for $(e_0,e_1,e_{\infty}) = (0, \frac12, \frac 1m)$ with $m \in \{3,4,6\}$,
and contains various improvements: It decomposes $f$ to obtain a smaller solution if possible, and compares exponent-differences to reduce the number of cases, see \cite{VJKThesis} for details.

\begin{appendix}
\section{Appendix}
\label{appendix}

We tabulate all $\text{Belyi}^{(1)}$ maps $O_1,\ldots,O_{20}$ with $|\,E_{_{\infty 24}}\,| = 5$,
and all $\text{Belyi}^{(1)}$ maps $P_1,\ldots,P_{12}$ with \linebreak $|\,E_{_{\infty 26}}\,| = 5$.
See \cite{dessindatabase} for explicit expressions in $\mathbb{Q}(s)(x)$ for each of these maps.

As one can see, there is only one orbit for each branching pattern, except for $(2^3)$, $(2^3)$, $(1^2,4)$ for which there is none. This means
 that there do not exist $g_0,g_1,g_t,g_{\infty} \in S_6$ for which each of $g_0,g_1$ is a product of 3 disjoint 2-cycles, $g_t$ is a 2-cycle, $g_{\infty}$ is a 4-cycle,
$g_0 g_1 g_t g_{\infty} = 1$, for which $<\hspace{-2pt}g_0,g_1,g_t,g_{\infty}\hspace{-2pt}>$ is transitive.

\begin{table}[H]
{\small
  \begin{tabular}{c|l|l|c|l||c|l | l|c|l}
 $n$ & branching pattern & \hspace{-2pt}name\hspace{-2pt}  & $|\cal{O}|$ & \hspace{-2pt}decomp.\hspace{-3pt} & $n$ & branching pattern& \hspace{-2pt}name\hspace{-2pt}  & $|\cal{O}|$ & \hspace{-2pt}decomp. \\
 \hline
   &                                   &                 &      &            &    &                                   &          &      &\\[-8pt]
2  & $(2)$, $(1^2)$, $(1^2)$           & $O_{1}$         & 1    & & 2  &     $(2)$, $(1^2)$, $(1^2)$           & $P_{1}$         & 1    & \\
   & $(1^2)$, $(1^2)$, $(2)$           & $O_{2}$         & 1    & &  &         $(1^2)$, $(1^2)$, $(2)$           & $P_{2}$         & 1    & \\
3  & $(3)$, $(1,2)$, $(1^3)$           & $O_{3}$         & 1    & &  3 &     $(3)$, $(1,2)$, $(1^3)$           & $P_{3}$         & 1    & \\
   & $(1,2)$, $(1,2)$, $(1,2)$         & $O_{4}$         & 4    &   &   &   $(1,2)$, $(1,2)$, $(1,2)$         & $P_{4}$         & 4    &\\
   & $(1^3)$, $(1,2)$, $(3)$           & $O_{5}$         & 1    &    &   &  $(1^3)$, $(1,2)$, $(3)$           & $P_{5}$         & 1    &\\
4  & $(1^2,2)$, $(1^2,2)$, $(4)$       & $O_{6}$         & 4    &  & 4 &     $(4)$, $(2^2)$, $(1^4)$           & $P_{6}$         & 1    &  $2 \circ 2$\\
   & $(4)$, $(2^2)$, $(1^4)$           & $O_{7}$         & 1    & $2 \circ 2$ & &   $(1,3)$, $(2^2)$, $(1^2,2)$       & $P_{7}$         & 3    &   \\
   & $(1,3)$, $(2^2)$, $(1^2,2)$       & $O_{8}$         & 3    &   &  &   $(2^2)$, $(2^2)$, $(1^2,2)$       & $P_{8}$         & 2    & $2 \circ 2$   \\
   & $(2^2)$, $(2^2)$, $(1^2,2)$       & $O_{9}$         & 2    & $2 \circ 2$ & &       $(1^2,2)$, $(2^2)$, $(1,3)$       & $P_{9}$         & 3    & \\
   & $(1^2,2)$, $(2^2)$, $(1,3)$       & $O_{10}$        & 3    &  &   &    $(1^2,2)$, $(2^2)$, $(2^2)$       & $P_{10}$        & 2    & $2 \circ 2$  \\
   & $(1^2,2)$, $(2^2)$, $(2^2)$       & $O_{11}$        & 2    & $2 \circ 2$ & &      $(1^4)$, $(2^2)$, $(4)$           & $P_{11}$        & 1    & $2 \circ 2$    \\
5  & $(1^2,3)$, $(1,2^2)$, $(1,4)$     & $O_{12}$        & 10   &  &   6  &  $(1^4,2)$, $(2^3)$, $(6)$         & $P_{12}$        & 2    & $2 \circ 3$  \\
   & $(1,2^2)$, $(1,2^2)$, $(1,4)$     & $O_{13}$        & 8    &  &   &   &&&\\
6  & $(1^2,4)$, $(2^3)$, $(1^2,4)$     & $O_{14}$        & 6    &  &   &    &&&\\
   & $(1,2,3)$, $(2^3)$, $(1^2,4)$     & $O_{15}$        & 9    &  &    &   &&&\\
   & $(2^3)$, $(2^3)$, $(1^2,4)$       & \hspace{2pt}---             & ---  &   &       &      &&&\\
   & $(1^3,3)$, $(2^3)$, $(2,4)$       & $O_{16}$        & 2    & $2 \circ 3$ & &     &&&\\
   & $(1^2,2^2)$, $(2^3)$, $(2,4)$     & $O_{17}$        & 4    & $2 \circ 3$  & &    &&&\\
8  & $(1^4,4)$, $(2^4)$, $(4^2)$       & $O_{18}$        & 2    & $2 \circ 2 \circ 2$  &&  &&&\\
   & $(1^3,2,3)$, $(2^4)$, $(4^2)$     & $O_{19}$        & 6    & $2 \circ 4$  &&     &&&\\
   & $(1^2,2^3)$, $(2^4)$, $(4^2)$     & $O_{20}$        & 4    & $2 \circ 2 \circ 2$ &&  &&&\\
 \end{tabular}
 }
\caption{$\text{Belyi}^{(1)}$ with $|\,E_{_{\infty 24}}\,| = 5$ resp. $|\,E_{_{\infty 26}}\,| = 5$} \label{table46}
\end{table}

\end{appendix}

\end{document}